\documentclass[a4paper,11pt]{article}%
\usepackage{amssymb,amsfonts,amsmath,amssymb,indentfirst,endnotes,rotating,arydshln,amsthm}
\usepackage[nohead]{geometry}
\usepackage[singlespacing]{setspace}
\usepackage[bottom]{footmisc}
\usepackage{graphicx,array,booktabs,bbm,multicol,colortbl,tabularx}
\usepackage{caption,subcaption}
\usepackage{authblk}
\usepackage{mathrsfs,mathtools}

\RequirePackage[OT1]{fontenc} \RequirePackage{amsthm,amsmath}
\RequirePackage[numbers]{natbib}
\RequirePackage[colorlinks,citecolor=blue,urlcolor=blue]{hyperref}
\RequirePackage{hypernat} \makeatletter
\newcommand{\field}[1]{\mathbb{#1}}

\newcommand{\X}{\field{X}}
\newcommand{\Y}{\field{Y}}
\newcommand{\E}{\field{E}}

\theoremstyle{example} \theoremstyle{remark} \theoremstyle{lemma}
\theoremstyle{definition} \theoremstyle{corol}
\theoremstyle{proposition} \theoremstyle{condition}
\theoremstyle{assumption}

\newtheorem{theorem}{\n{Theorem}}[section]

\newtheorem{remark}{\n{Remark}}[section]
\newtheorem{lemma}{\n{Lemma}}[section]

\newtheorem{corollary}{\n{Corollary}}[section]

\newcommand{\Rmnum}[1]{\expandafter\romannumeral #1}

\font\n=cmcsc12

 
\newcommand{\dotp}[2]{\left\langle #1, #2\right\rangle}
\newcommand{\argmin}{\mathop{\rm argmin~}}

\newcommand{\ve}{\varepsilon}

\newcommand{\wt}{\widetilde}
\newcommand{\wht}{\widehat}

\newcommand{\abs}[1]{\left\vert#1\right\vert}

\newcommand{\norm}[1]{\left\Vert#1\right\Vert}

\newcommand \bbP{\mathbb{P}}
\newcommand \bbE{\mathbb{E}}

\def\m{\mathcal}
\def\mb{\mathbb}
\def\mr{\mathrm}

\newcommand{\ind}{\mathbbm{1}}

\def\T{{ \mathrm{\scriptscriptstyle T} }}
\def\So{{ \mathrm{\scriptscriptstyle S} }}
\def\H{{ \mathrm{\scriptscriptstyle H} }}

\newcommand{\bmp}{\begin{minipage}}
\newcommand{\emp}{\end{minipage}}
\newcommand{\ben}{\begin{eqnarray*}}
\newcommand{\non}{\end{eqnarray*}}
\newcommand{\1}{\\[1ex]}

\newcommand{\bc}{\begin{center}}
\newcommand{\ec}{\end{center}}

\begin{document}
\title{\Large Frequentist coverage and sup-norm convergence rate in Gaussian process regression}
\author[1]{Yun Yang \thanks{yyang@stat.fsu.edu}}
\author[2]{Anirban Bhattacharya\thanks{anirbanb@stat.tamu.edu}}
\author[2]{Debdeep Pati\thanks{debdeep@stat.tamu.edu}}
\affil[1]{Department of Statistics, Florida State University}
\affil[2]{Department of Statistics, Texas A\&M University}
\date{\vspace{-2em}}

\maketitle

\begin{abstract}
Gaussian process (GP) regression is a powerful interpolation technique due to its flexibility in capturing non-linearity. In this paper, we provide a general framework for understanding the frequentist coverage of point-wise and simultaneous Bayesian credible sets in GP regression.  As an intermediate result, we develop a Bernstein von-Mises type result under supremum norm in random design GP regression.   
Identifying both the mean and covariance function of the posterior distribution of the Gaussian process as regularized $M$-estimators, we show that the sampling distribution of the posterior mean function and the centered posterior distribution can be respectively approximated by two population level GPs. By developing a {\em comparison inequality} between two GPs, we provide exact characterization of frequentist coverage probabilities of Bayesian point-wise credible intervals and simultaneous credible bands of the regression function. Our results show that inference based on GP regression tends to be conservative; when the prior is under-smoothed, the resulting credible intervals and bands have minimax-optimal sizes, with their frequentist coverage converging to a non-degenerate value between their nominal level and one.  As a byproduct of our theory, we show that the GP regression also yields minimax-optimal posterior contraction rate relative to the supremum norm, which provides a positive evidence to the long standing problem on optimal supremum norm contraction rate in GP regression. 
 
\end{abstract}


\section{Introduction and preliminaries}
Gaussian process (GP) regression is a popular Bayesian procedure for learning an infinite-dimensional function $f\in\m F$ in the nonparametric regression model
\begin{align*}
y=f(x) + \varepsilon,\quad \varepsilon\sim\m N(0,\sigma^2),
\end{align*}
where $x\in\m X\in [0,1]^d$ is the covariate variable and $y\in\mb R^d$ the response variable. Through specifying a GP as the prior measure over the infinite dimensional parameter space $\m F$, Bayes†rule yields a posterior measure for $f$ that can be used to either construct a point estimator as the posterior mean, or characterize estimation uncertainties through the corresponding point-wise credible intervals and simultaneous credible bands. 
Examples of the wide usages of GP include computer experiment emulations \citep{sacks1989design,currin1991bayesian,kennedy2001bayesian}, spatial data modeling \citep{cressie2015statistics,banerjee2008gaussian}, geostatistical kriging \citep{stein2012interpolation,matheron1973intrinsic}, and machine learning \citep{rasmussen2006gaussian}. 

Despite its long-standing popularity, formal investigations into theoretical properties of GP regression from a frequentist perspective assuming the existence of a true data-generating function $f^\ast$ is a much more recent activity. 
A majority of this recent line of work has been directed towards understanding large sample performance of point estimation in the form of proving posterior consistency and deriving posterior contraction rates. As an incomplete survey, 
\citep{van2008rates,vaart2011information} provide a general framework for studying the rate of posterior contraction of GP regression, and derive the contraction rates for several commonly used covariance kernels. In \cite{vandervaart2009} and \cite{bhattacharya2014anisotropic}, the authors show that by putting hyper priors over inverse bandwidth parameters in the square-exponential covariance kernel, GP regression can adapt to the unknown (anisotropic) smoothness levels of the target function. \citep{yang2016bayesian} show that a class of GP priors with Euclidean-metric based covariance kernels can additionally adapt to unknown low-dimensional manifold structure even with a moderate-to-high dimensional ambient space $\m X$ in which the covariate lives. 

There is a comparatively limited literature on the validity of conducting statistical inference in GP regression, or more generally, in Bayesian nonparametric procedures, from a frequentist perspective. Uncertainty quantification for GPs plays an important role---even more important than point estimation itself---in some applications such as design of experiments \citep{burnaev2015adaptive} and risk management \cite{scansaroli2012stochastic}, and hence it is of interest to investigate the frequentist validity of such. However, unlike finite dimensional parametric models, where a Bernstein von Mises (BvM) theorem under the total variation metric holds fairly generally and guarantees asymptotically nominal coverage of Bayesian credible intervals, the story of the frequentist validity of credible intervals/bands in infinite dimensional models is more complicated and delicate \cite{cox1993analysis,freedman1999wald,johnstone2010high,leahu2011bernstein}. For the Gaussian white noise model, \citep{castillo2013nonparametric} showed that a Bernstein von Mises (BvM) theorem holds in weaker topologies for some natural priors, which yields the correct asymptotic coverage of credible sets based on the weaker topology; however their result is not applicable for understanding the coverage of credible intervals/bands. A similar result for the stronger $L^\infty$-norm can be found in \cite{castillo2014bernstein}.  In the context of linear inverse problems with Gaussian priors, \cite{knapik2011bayesian} showed that asymptotic coverage of credible sets can be guaranteed by under-smoothing the prior compared to the truth.  \cite{szabo2015frequentist} investigated the frequentist coverage of Bayesian credible sets in the Gaussian white noise model, and showed that depending on whether the smoothness parameter in the prior distribution matches that in the truth, the coverage of the corresponding credible sets can either be significantly below its nominal confidence level, or converge to one as the sample size increases, even though the nominal level is fixed at a constant value. \cite{szabo2015frequentist} also investigated the adaptivity of the credible sets to unknown smoothness using an empirical Bayes approach; see also \cite{serra2017adaptive,ray2014adaptive}. 

A majority of the work discussed in the previous paragraph focuses on the Gaussian white noise model or its equivalent infinite dimensional normal mean formulation, where appropriate conjugate Gaussian priors lead to an analytically tractable posterior. 
A major difficulty of calculating the nominal coverage probabilities of Bayesian point-wise credible intervals and simultaneous credible bands in GP regression lies in the fact that the covariance structure in the posterior GP of $f$ involves unwieldy matrix inverses that are complicated for direct analysis. Moreover, when the design is random, the randomness in the covariance structure further complicates the analysis. A relevant work in this regard is \cite{yoo2016supremum}, where the authors derive the posterior contraction rates under the supremum norm for GP regression induced by tensor products of B-splines, and identify conditions under which the coverage of the simultaneous credible band tends to one as the sample size tends to infinity. However, their result requires the nominal level of the credible bands to also tend to one, and is based on a key property of their prior distribution that the resulting posterior covariance function can be sandwiched by two identity matrices with similar scalings, preventing it to be applicable to a broader class of GP priors. Also relevant to the present discussion is \cite{sniekers2015credible}, who obtained  similar results for point-wise credible intervals using a scaled Brownian motion prior. The authors exploit explicit formulas for the eigenvalues and eigenfunctions of the covariance kernel of a Brownian motion when the design points are on a regular grid, which along with other properties specific to Brownian motion can be used to linearize the posterior mean and variance, aiding a direct analysis. 

The goal of this article is to provide a general framework for understanding the frequentist coverage of Bayesian credible intervals and bands in GP regression by proving a BvM type theorem under the supremum norm in random design nonparametric regression. Towards this goal, we first show that the sampling distribution of the posterior mean function can be well-approximated by a GP over the covariate space $\m X$ with an explicit expression for its covariance function. Second, we find a tractable population level GP approximation to the centered posterior GP whose covariance function is non-random and also admits an explicit expression. The frequentist coverage of the Bayesian credible intervals and bands are derived from an interplay between these two population level GPs. 
A salient feature of our technique is that it provides an explicit expression of the coverage along with explicit finite sample error bounds---it is non-asymptotic and applies to any nominal level. Interestingly, we find that when the prior is under-smoothed, the Bayesian credible intervals and bands are always moderately conservative in the sense that its frequentist coverage is always higher than its nominal level, and converges to a {\em fixed number} (with explicit expression) between its nominal level and one as the sample size grows to infinity. For example, when the covariate is one-dimensional and uniformly distributed over the unit interval, and the prior smoothness level is $2(3)$, the asymptotic frequentist coverage of a $95\%$ credible interval is 0.976(0.969). This phenomenon is radically different from existing results, where the frequentist coverage is either tending to zero or one but never converges to a non-degenerate value. As a byproduct of our theory, we show that the GP regression also yields minimax-optimal posterior contraction rate relative to the supremum norm, which provides positive evidence to the long standing problem about supremum norm contraction rate of general Bayesian nonparametric procedures. 

In our proofs, we employ the equivalent kernel representation \citep{rasmussen2006gaussian,silverman1984spline} of the kernel ridge regression estimator to establish a novel local linear expansion of the posterior mean function relative to the supremum norm, which is of independent interest and builds a link between GP regression and frequentist kernel type estimators \cite{gine2008uniform}. This local linear expansion can be utilized to show the limiting GP approximation to the sampling distribution of the posterior mean function. Towards the proof of approximating the posterior GP with a population level GP, we develop a new Gaussian comparison inequality that provides explicit bounds on the Kolmogorov distance between two GPs through the supremum norm difference between their respective covariance functions.

Overall, our results reveals delicate interplay between frequentist coverage of Bayesian credible sets and proper characteristics of the prior measure in infinite-dimensional Bayesian procedures. We validate GP regression for conducting statistical inference by showing that as long as the prior measure is not over-smoothed, Bayesian credible sets always provide moderately conservative uncertainty quantification with minimax-optimal sizes. 

We begin the technical development by introducing notation used throughout the paper. A summary of all the major notations are provided in Table \ref{tab:notes} in \S \ref{sec:proof} for the reader's convenience.

\subsection{Notation}


Let $\mathbb{H}_1,  \mathbb{H}_2$ be normed linear spaces.  The Fr\'{e}chet derivative of an operator $A$ at the point 
$a \in \mathbb{H}_1$ is the  bounded linear operator denoted $DA (a):  \mathbb{H}_1 \to  \mathbb{H}_2$ which satisfies 
\begin{eqnarray}\label{eq:fdef}
\lim_{h \to 0} \frac{\norm{A(a + h) -  A(a) -  DA(a) h}}{\norm{h}} = 0.  
\end{eqnarray}
In particular, when $\mb H_1 = \mb R^n, \mb H_2 = \mb R^m$, the Fr\'{e}chet derivative $DA(a)$  is the Jacobian of $A$, 
a linear operator which is represented by an $m \times n$ matrix $(\partial/\partial x_j A_i)$. 
 
Throughout $C, C', C_1, C_2, \ldots$ are generically used to denote positive constants whose values might change from one line to another, but are independent from everything else.  We use $\lesssim$ and $\gtrsim$ denote inequalities up to a constant multiple; we write $a \asymp b$  when both $a \lesssim b$ and $a \gtrsim b$ hold. For $\alpha > 0$, let $\lfloor \alpha \rfloor$ denote the largest integer strictly smaller than $\alpha$.

\subsection{Review of RKHS}

We recall some key facts related to reproducing kernel Hilbert spaces (RKHS); further details and proofs can be found in Chapter 1 of \cite{wahba1990spline}. 
\noindent Let $\m X $ denote a general index set. A symmetric function $K : \m X \times \m X \to \mb R$ is {\em positive definite} (p.d.) if for any $n \in \mb N$, $a_1, \ldots, a_n \in \mb R$, and $t_1, \ldots, t_n \in \m X$, 
\begin{align*}
\sum_{i=1}^n \sum_{j=1}^n a_i a_j K(t_i, t_j) > 0. 
\end{align*}
A (real) RKHS $\m H$ is a Hilbert space of real-valued functions on $\m X$ such that for any $t \in \m X$, the evaluation function $L_t : \m H \to \mb R$; $f \mapsto f(t)$, is a bounded linear functional, i.e., there exists a constant $M_t > 0$ such that 
$$
|L_t f| = |f(t)| \le M_t \norm{f}_{\m H}, \quad f \in \m H. 
$$
In the above display, $\norm{f}_{\m H} = \sqrt{ \dotp{f}{f}_{\m H} }$ is the Hilbert space norm. Since the evaluation maps are bounded linear, it follows from the Riesz representation theorem that for each $t \in \m X$, there exists an element $K_t \in \m H$ such that $f(t) = \dotp{f}{K_t}_{\m H}$. $K_t$ is called the {\em representer of evaluation} at $t$, and the kernel $K(s, t) = K_t(s)$ is easily shown to be p.d. 
Conversely, given a p.d. function $K$ on $\m X \times \m X$, one can construct a unique RKHS on $\m X$ with $K$ as its reproducing kernel. Given a kernel $K$, we shall henceforth let $K_t(\cdot) = K(\cdot, t)$. 

Let $L^2(\m X)$ denote the space of square integrable functions $f:\m X \to \mb R$ with $\int_{\m X} f^2(x) dx < \infty$. We denote $\dotp{f}{g}_{L^2(\m X)} := \int_{\m X} f(x) g(x) dx$ the usual inner product on $\m X$. If a p.d. kernel $K(s, t)$ is continuous and $\int_{\m X} \int _{\m X} K^2(s, t) ds dt < \infty$, then by Mercer's theorem, there exists an orthonormal sequence of continuous eigenfunctions $\{ \psi_j \}_{j=1}^{\infty}$ in $L^2(\m X)$ with eigenvalues $\mu_1 \ge \mu_2 \ldots \ge 0$, and 
\begin{align*}
& \int K(s, t) \psi_j(s) ds  = \mu_j \psi_j(t), \quad j = 1, 2, \ldots, \\
& K(s, t) = \sum_{j=1}^{\infty} \mu_j \psi_j(s) \psi_j(t), \quad \int_{\m X} \int _{\m X} K^2(s, t) ds dt = \sum_{j=1}^{\infty} \mu_j^2 < \infty. 
\end{align*}
The RKHS $\m H$ determined by $K$ consists of functions $f \in L^2(\m X)$ satisfying $\sum_{j=1}^{\infty} f_j^2/\mu_j < \infty$, where $f_j = \dotp{f}{\psi_j}_{L^2(\m X)}$. Further, for $f, g \in \m H$, 
\begin{align}\label{eq:rkhs_series}
\norm{f}_{\m H}^2 = \sum_{j=1}^{\infty} \frac{f_j^2}{\mu_j}, \quad \dotp{f}{g}_{\m H} = \sum_{j=1}^{\infty} \frac{f_j g_j}{\mu_j}. 
\end{align}
For stationary kernel $K$ with $\m X=[0,1]$ and $\mu$ being the Lebesgue measure, we can always choose $\psi_{2j-1}(s)=\sin(\pi j s)$ and $\psi_{2j}(s)=\cos(\pi j s)$, $j=1,\,2,\ldots$, by expanding $K(s,\,t)$ into a Fourier series over $[-1,\,1]$ and exploiting the identity $\cos(x-y)=\cos(x)\cos(y)+\sin(x)\sin(y)$. Under these choices, we also have $\mu_{2j-1}=\mu_{2j}$ for $j=1,\,2,\ldots$ (more details can be found in Appendix A of \cite{yang2017bayesian}).

\section{Framework}
To begin with, we introduce Gaussian process regression and draw its connection (for both the posterior mean function and posterior covariance function) with kernel ridge regression (KRR) in \S \ref{sec:gp_reg}. In \S \ref{subsec:EK}, we introduce the key mathematical tool in our proofs---equivalent kernel representation of the kernel ridge regression. In \S \ref{sec:supnorm_krr}, we present our first result on the local linear expansion of the KRR estimator relative to the supremum norm, indicating the asymptotic equivalence between the KRR estimator with a carefully constructed kernel type estimator.

\subsection{GP regression}\label{sec:gp_reg}
For easy presentation, we focus on the univariate regression problem where $\m X\subset \mb  R$, and our results can be straightforwardly extend to multivariate cases.
Let $(Y_i, X_i), i = 1, \ldots, n$, be i.i.d., with $X_i \in \m X$ and $Y_i \in \mb R$, with joint density $\rho_{Y, X}(y, x) \propto \phi_\sigma(y - f^\ast(x)) \ \ind_{\m X}(x)$, where $\phi_\sigma$ denotes the $\m N(0, \sigma^2)$ density, and $f^\ast : \m X \to \mb R$ is the unknown function of interest. Our goal is to estimate and perform inference on $f^\ast$ based on the data $\mb{D}_n =  \{(Y_i, X_i), i = 1, \ldots, n\}$. 
We assume $\sigma^2$ is known throughout the paper. 

We consider a nonparametric regression model 
\begin{align}\label{eq:npreg}
Y_i = f(X_i) + \varepsilon_i, \quad \varepsilon_i \sim N(0, \sigma^2), 
\end{align}
and assume a mean-zero GP prior on the regression function $f$, $f \sim \mbox{GP}(0, \sigma^2 (n \lambda)^{-1} K)$, where $K$ is a positive definite kernel and $\lambda > 0$ is a tuning parameter to be chosen later. 

By conjugacy, it is easy to check that the posterior distribution of $f$ is also a GP, $f \mid \mb D_n \sim \mbox{GP}(\widehat{f}_n, \wt{C}_n^B)$, with mean function $\widehat{f}_n(\cdot)$ and covariance function $\wt{C}_n^B(\cdot, \cdot)$,
\begin{align}\label{eq:meancov}
\widehat{f}_n(x) =  \bbE[f(x) \mid \mb{D}_n], \quad \wt{C}_n^B(x,x')=\mbox{Cov}[f(x),\,f(x')\mid \mb D_n], \quad x, x' \in \m X. 
\end{align}
Since the posterior distribution is a GP, it is completely characterized by $\widehat{f}_n$ and $\wt{C}_n^B$. 
These quantities can be explicitly calculated; we introduce some notation to express these succinctly. For vectors $\E = (e_1, \ldots, e_m)^{\T}$ and $\E' = (e_1', \ldots, e_r')^{\T}$, let $K(\E, \E')$ denote the $m \times r$ matrix $(K(e_i, e_j'))$. Also, let $\X = (X_1, \ldots, X_n)^{\T}$ and $\Y = (Y_1, \ldots, Y_n)^{\T}$. With these notations, 
\begin{align}
\widehat{f}_n(x)  & =  K(x, \X)\, [K(\X, \X) + n \lambda \, \mr I_n]^{-1} \Y,  \label{eq:postme}\\
\wt{C}_n^B(x,x') & = \sigma^2\, (n\lambda)^{-1}\bigg\{ K(x, x') -  K(x, \X)\, [K(\X, \X) + n\lambda \, \mr I_n]^{-1} K(\X, x') \bigg\} \label{eq:postcov}. 
\end{align}
In particular, the posterior variance function is 
$$\mbox{Var}[f(x) \mid \mathbb{D}_n] = \sigma^2\,(n\lambda)^{-1}\bigg\{K(x, x) -  K(x, \X)\, [K(\X, \X) + n\lambda\, \mr I_n]^{-1} K(\X, x)\bigg\}.$$ 

The presence of the inverse kernel matrix in equations \eqref{eq:postme} and \eqref{eq:postcov} renders analysis of the GP posterior unwieldy. A contribution of this article is to recognize \emph{both} the mean function and covariance function as regularized $M$-estimators and use a equivalent kernel trick to linearize the solutions that avoids having to deal with matrix inverses. 
It is well-known (see, e.g. Chapter 6 of \cite{rasmussen2006gaussian}) that the posterior mean $\wht{f}_n$ coincides with the kernel ridge regression (KRR) estimator
\begin{align}\label{eq:krr}
\wht{f}_{n, \lambda} = \argmin_{f \in \m H} \ell_{n, \lambda}(f), \quad \ell_{n, \lambda}(f) := \bigg[ \frac{1}{n} \sum_{i=1}^n (Y_i - f(X_i))^2 + \lambda \norm{f}_{\m H}^2 \bigg], 
\end{align}
when the RKHS $\m H$ corresponds to the prior covariance kernel $K$. The objective function $\ell_{n, \lambda}$ in \eqref{eq:krr} combines the average squared-error loss with a squared RKHS norm penalty weighted by the prior precision parameter $\lambda$. It follows from the representer theorem for RKHSs \cite{wahba1990spline} that the solution to \eqref{eq:krr} belongs to the linear span of the kernel functions $\{ K(\cdot, X_i) \}_{i=1}^n$. Given this fact, solving \eqref{eq:krr} only amounts to solving a quadratic program, and the solution coincides with $\wht{f}_n$ in \eqref{eq:postme}. 

A novel observation aiding our subsequent development is that the posterior covariance function $\wt{C}_n^B$ can be related to the bias of a noise-free KRR estimator. Write 
$\sigma^{-2}\,n\lambda \, \wt{C}^B_n(x,x')$ as
\begin{align*}
\sigma^{-2}\,n \lambda\, \wt{C}^B_n(x,x')  =  K_{x}(x') - \wht{K}_{x}(x'),
\end{align*}
where $K_x(\cdot) = K(x, \cdot)$ as defined earlier, and $\wht{K}_x(\cdot) = K(\cdot, \X)\, [K(\X, \X) + n\lambda \, \mr I_n]^{-1} K(\X, x)$. Comparing with \eqref{eq:postme} and \eqref{eq:krr}, it becomes apparent that $\wht{K}_x$ is the solution to the following KRR problem with noiseless observations of $K_x$ at the random design points $\{X_i\}_{i=1}^n$, 
\begin{align}\label{eq:var_krr}
\widehat{K}_{x} = \argmin_{g \in \m H} \bigg[ \frac{1}{n} \sum_{i=1}^n (Z^x_i - g(X_i))^2+ \lambda \norm{g}_{\m H}^2 \bigg],
\end{align}
where $Z^x_i = K(x, X_i) = K_x(X_i)$. 

To summarize, the posterior mean corresponds to the usual KRR estimator, and an appropriately scaled version of the posterior covariance function can be recognized as the bias of a noiseless KRR estimator. This motivates us to study the performance of KRR estimators in the supremum norm, which to best of our knowledge, hasn't been carried out before. For past work on risk bounds for the KRR estimator in other norms, refer to \cite{mendelson2002geometric,zhang2005learning,steinwart2009optimal,dicker2017kernel,yang2017non}. 

\subsection{Equivalent kernel representation of the KRR estimator} \label{subsec:EK}

We first introduce an equivalent-kernel formulation that allows us to linearize the bias of a KRR estimate. 
Let $\m H \subset L^2(\m X)$ be an RKHS, with reproducing kernel $K$ having eigenfunctions $\{ \psi_j \}$ and eigenvalues $\{\mu_j\}$ with respect to $L^2(\m X)$. Fix $\lambda > 0$. Define a new inner product on $\m H$ as 
\begin{align}\label{eq:newdot_def}
\dotp{f}{g}_{\lambda} := \dotp{f}{g}_{L^2(\m X)} + \lambda \dotp{f}{g}_{\m H},
\end{align}
and let $\norm{f}_{\lambda} = \sqrt{ \dotp{f}{f}_{\lambda}}$. Observe that $(\m H, \dotp{\cdot}{\cdot}_{\lambda})$ is again an RKHS, since for any $x \in \m X$, $|f(x)| \le C \norm{f}_{\m H} \le C' \norm{f}_{\lambda}$ for some constant $C' > 0$ depending on $C$ and $\lambda$, proving the boundedness of the evaluation maps. 

Let $f = \sum_{j=1}^{\infty} f_j \psi_j$ and $g = \sum_{j=1}^{\infty} g_j \psi_j$ be elements of $L^2(\m X)$. Then, 
\begin{align*}
\dotp{f}{g}_{L^2(\m X)} + \lambda \dotp{f}{g}_{\m H} = \sum_{j=1}^{\infty} f_j g_j +\lambda \sum_{j=1}^{\infty} \frac{f_j g_j}{\mu_j} = \sum_{j=1}^{\infty} \frac{f_j g_j}{\nu_j}, 
\end{align*}
where 
\begin{align}\label{eq:def_nu}
\nu_j = \frac{1}{1 + \frac{\lambda}{\mu_j}} = \frac{\mu_j}{\lambda + \mu_j}, \quad j = 1, 2, \ldots
\end{align}
Hence, $(\m H, \dotp{\cdot}{\cdot}_{\lambda})$ consists of  
$\big\{ f = \sum_{j=1}^{\infty} f_j \psi_j \in L^2(\m X): \sum_{j=1}^{\infty} f_j^2/\nu_j < \infty \big\}$, 
with 
\begin{align}\label{eq:dotp_nu}
\dotp{f}{g}_{\lambda} = \sum_{j=1}^{\infty} \frac{f_j g_j}{\nu_j}. 
\end{align}
The reproducing kernel associated with the new RKHS is thus 
\begin{align}\label{eq:new_kernel}
\wt{K}(s, t) := \sum_{j=1}^{\infty} \nu_j \psi_j(s) \psi_j(t). 
\end{align}
As before, we let $\wt{K}_s(\cdot)$ denote the representer of the evaluation map, i.e., $\wt{K}_s(\cdot) = \wt{K}(s, \cdot)$, so that for any $g \in \m H$, $g(s) = \langle{g, \,}{\wt{K}_s}\rangle_{\lambda}$. The kernel $\wt{K}$ is known as the {\em equivalent kernel} (see, e.g., Chapter 7 of \cite{rasmussen2006gaussian}), motivated by the notion of equivalent kernels in the spline smoothing literature \cite{silverman1984spline}. We comment more on connections with the equivalent kernel literature once we represent the KRR problem in terms of the equivalent kernel.

Before proceeding further, we introduce two operators which are routinely used subsequently. First, define a linear operator $F_{\lambda}: \m H \to \m H$ given by 
\begin{align*}
F_{\lambda} g(t) = \int g(s) \wt{K}(s, t) ds.
\end{align*}
We shall often use the abbreviation $F_{\lambda} g = \int g(s) \wt{K}_s ds$. The operator $F_{\lambda}$ is easily recognized as a convolution with the equivalent kernel $\wt{K}$. If $g = \sum_{j=1}^{\infty} g_j \psi_j$, a straightforward calculation yields $F_{\lambda} g = \sum_{j=1}^{\infty} \nu_j g_j \psi_j$. Thus, for functions $f, g \in \m H$, it follows from \eqref{eq:dotp_nu} that 
\begin{align*}
\dotp{f}{F_{\lambda} g}_{\lambda} = \sum_{j=1}^{\infty} f_j g_j = \dotp{f}{g}_{L^2(\m X)}. 
\end{align*}
The above display also immediately tells us that $F_{\lambda}$ is a self-adjoint operator, i.e., $\dotp{f}{F_{\lambda} g}_{\lambda} = \dotp{F_{\lambda}f}{g}_{\lambda}$ for all $f, g \in \m H$. Define another self-adjoint operator $P_{\lambda}: \m H \to \m H$ given by $P_{\lambda} = \mbox{id} - F_{\lambda}$, where $\mbox{id}$ is the identity operator on $\m H$. Then, it follows from the previous display and \eqref{eq:newdot_def} that 
\begin{align}\label{eq:Plam_id}
\dotp{f}{P_{\lambda} g}_{\lambda} = \lambda \dotp{f}{g}_{\m H}, \quad f, g \in \m H. 
\end{align}
Having developed the necessary groundwork, let us turn our attention back to the KRR estimate. Recall the objective function $\ell_{n, \lambda}$ defined in \eqref{eq:krr}. Writing $f(X_i) = \langle f, \wt{K}_{X_i} \rangle_{\lambda}$, and $\norm{f}_{\m H}^2 = \dotp{f}{P_{\lambda} f}_{\lambda}$ using \eqref{eq:Plam_id}, we can express 
\begin{align*}
\ell_{n, \lambda}(f) = \bigg[ \frac{1}{n} \sum_{i=1}^n (Y_i - \langle f, \wt{K}_{X_i} \rangle_{\lambda})^2+ \dotp{f}{P_{\lambda}f}_{\lambda} \bigg]. 
\end{align*}
Viewing $\ell_{n, \lambda} : (\m H, \dotp{\cdot}{\cdot}_{\lambda}) \to \mb R$ and performing a Fr\'{e}chet differentiation with respect to $f$, one obtains a score equation for the KRR estimate $\wht{f}_{n, \lambda}$ as, 
$$
S_{n, \lambda}( \wht{f}_{n, \lambda}) = 0, 
$$
where $S_{n, \lambda} : \m H \to \m H$ is given by 
\begin{align}\label{eq:snlam}
S_{n, \lambda}(f) = \frac{1}{n} \sum_{i=1}^n (Y_i - f(X_i)) \wt{K}_{X_i} - P_{\lambda} f. 
\end{align}
Define $S_{\lambda}(f) :=  \bbE_{f^\ast} [ S_{n, \lambda}(f)]$ to be the population version of the score equation, where the expectation is assumed with respect to the true joint density $\rho_{Y, X}$. Recall the convolution operator $F_{\lambda}$. We then have,
\begin{align}\label{eq:slam}
S_{\lambda}(f)  = \int_{\m X} \{f^\ast(x) -  f(x)\} \,\wt{K}_x \,dx -  P_{\lambda} f = F_{\lambda}(f^\ast - f) - P_{\lambda} f = F_{\lambda} f^\ast - f, 
\end{align}
since $F_{\lambda} f + P_{\lambda} f  = f$ by definition. It is immediate that $S_{\lambda}(F_{\lambda} f^\ast) = 0$, which therefore implies that the function $F_{\lambda} f^\ast$ is the solution to the population level score equation. We shall henceforth refer to $F_{\lambda} f^\ast$ as the population-level KRR estimator. 

In the above treatment, we first differentiated the objective function and then took an expectation with respect to the true distribution $\rho_{Y,X}$ to arrive at the population-level KRR estimator $F_{\lambda} f^\ast$. One arrives at an identical conclusion if the expectation of the objective function $\ell_{n, \lambda}$ is minimized, which is equivalent to minimizing 
$$
\norm{f^\ast - f}_{L^2(\m X)}^2 + \lambda \norm{f}_{\m H}^2 = \sum_{j=1}^{\infty} (f^\ast_j - f_j)^2 + \lambda \sum_{j=1}^{\infty} f_j^2/\lambda_j.
$$
Solving for $f_j$, one obtains $f_j = \nu_j f_j^*$, and hence $f = \sum_{j=1}^{\infty} \nu_j f_j^* \psi_j = F_{\lambda} f^*$. This approach thus also  leads to the equivalent kernel; see Chapter 7 of \cite{rasmussen2006gaussian} for a detailed exposition along these lines.

\subsection{Sup-norm bounds for the KRR estimator}\label{sec:supnorm_krr}

We now use the equivalent kernel representation to derive error bounds in the supremum norm between a KRR estimator and its target function. We first lay down two kinds of parameter space considered for the true function $f^*$. Recall the orthonormal basis $\{\psi_j\}_{j=1}^{\infty}$ of $L^2(\m X)$. 

For any $\alpha > 1/2$ and $B > 0$, define 
\begin{align}\label{eq:Sob}
\Theta_{\So}^{\alpha}(B) = \big\{ f = \sum_{j=1}^{\infty} f_j \psi_j \in L^2(\m X) : \sum_{j=1}^{\infty} j^{2 \alpha} f_j^2 \le B^2 \big\}.  
\end{align}
For integer $\alpha$ and the Fourier basis, $\Theta_{\So}^{\alpha}$ corresponds to the $\alpha$-smooth Sobolev functions with absolutely continuous $(\alpha - 1)$ derivatives and whose $\alpha$th derivative has uniformly bounded $L^2$ norm. 

Next, for any $\alpha > 0$ and $B > 0$, define 
\begin{align}\label{eq:Hol}
\Theta_{\H}^{\alpha}(B) = \big\{ f = \sum_{j=1}^{\infty} f_j \psi_j \in L^2(\m X) : \sum_{j=1}^{\infty} j^{\alpha} |f_j| \le B \big\}.  
\end{align}
The $\H$-subscript is used to indicate a correspondence of this class of functions with $\alpha$-H\"{o}lder functions. Indeed, under the Fourier basis, if $f \in \Theta_{\H}^{\alpha}$, then $f$ has $\lfloor \alpha \rfloor$ continuous derivatives, and 
$$
\abs{ f^{(\lfloor \alpha \rfloor)}(x) - f^{(\lfloor \alpha \rfloor)}(x') } \lesssim \sum_{j=1}^{\infty} j^{\lfloor \alpha \rfloor} |f_j| \ j^{\alpha - \lfloor \alpha \rfloor} |x - x'|^{\alpha - \lfloor \alpha \rfloor} \lesssim |x - x'|^{\alpha - \lfloor \alpha \rfloor},
$$
which implies that the $\lfloor \alpha \rfloor$th derivative of $f$ is Lipschitz continuous of order $\alpha - \lfloor \alpha \rfloor$. \\

We next lay down some standard technical conditions on the eigenfunctions and eigenvalues of the kernel $K$. 
\paragraph{Assumption (B):} There exists global constants $C_{\psi}, L_{\psi} > 0$ such that the eigenfunctions $\{\psi_j\}_{j=1}^{\infty}$ of $K$ satisfy $|\psi_j(t)| \le C_{\psi}$ for all $j \ge 1, \, t \in \m X$, and $|\psi_j(t) - \psi_j(s)| \le L_{\psi} \, j |t - s|$ for all $j \ge 1$ and $t, s \in \m X$. 

\paragraph{Assumption (E):} The eigenvalues $\{\mu_j\}_{j=1}^{\infty}$ of $K$ satisfy $\mu_j \asymp j^{-2 \alpha}$ for some $\alpha > 1/2$. \\[1ex]

As a motivating example, the Mat\'{e}rn kernel with smoothness index $(\alpha - 1/2)$ satisfies {\bf (B)} and {\bf (E)} when expanded with respect to the Fourier basis 
$$\psi_{2j-1}(x) = \sin(\pi j x), \ \psi_{2j}(x) = \cos(\pi j x), \ j = 1, 2, \ldots$$
on $\m X = [0, 1]$; see \cite{bhattacharya2015wass,yang2017bayesian} for more details.  Observe that the RKHS associated with a Gaussian process associated with Mat\'{e}rn kernel with smoothness index $(\alpha - 1/2)$ (in the multivariate case, $\alpha-d/2$, with $d$ the dimensionality of $\m X$) is  $\Theta_{\So}^{\alpha}$.  As a passing comment, this space does not contain functions with smoothness less than $\alpha$, which includes functions with smoothness $\alpha -1/2$. Although for concreteness we focus on kernels with polynomially decaying eigenvalues in the paper, our theory is also applicable to other kernel classes, such as squared exponential kernels and polynomial kernels.

Observe that $\nu_j = \mu_j/(\mu_j + \lambda) \asymp 1/(1 + \lambda j^{2 \alpha})$ under Assumption (E). Bounding the sums by integrals, the following facts are easily observed and used repeatedly in the sequel, 
\begin{align}\label{eq:nusum}
\sum_{j=1}^{\infty} \frac{1}{1 + \lambda j^{2 \alpha}} \lesssim \lambda^{-1/(2 \alpha)}, \quad \sum_{j=1}^{\infty} \frac{1}{(1 + \lambda j^{2 \alpha})^2} \lesssim \lambda^{-1/(2 \alpha)}. 
\end{align}

Recall the KRR estimator $\wht{f}_{n, \lambda}$ from \eqref{eq:krr}. Let $w_i = Y_i - f^\ast(X_i)$, so that $w_i \sim \mbox{N}(0, \sigma^2)$ are independent for $i = 1, \ldots, n$ and also independent of $(X_1, \ldots, X_n)^{\T}$. Recall the operators $F_{\lambda}$ and $P_{\lambda}$ from \S \ref{subsec:EK}. For $f = \sum_{j=1}^{\infty} f_j \psi_j \in L^2(\m X)$, we continue to denote 
\begin{align}\label{eq:PFL}
F_{\lambda} f = \sum_{j=1}^{\infty} \nu_j f_j \, \psi_j, \quad P_{\lambda} f = \sum_{j=1}^{\infty} (1 - \nu_j) f_j \, \psi_j.
\end{align}
Since $\nu_j \in (0, 1)$ for all $j \ge 1$, $F_{\lambda} f$ and $P_{\lambda} f$ are elements of $L^2(\m X)$ for any $f \in L^2(\m X)$. 

We now state a theorem which bounds the sup-norm distance between $\wht{f}_{n, \lambda}$ and the true function $f^\ast$. To state the theorem in its most general form, we don't make any smoothness assumptions of $f^\ast$ yet and state a high probability bound on $\| \wht{f}_{n, \lambda} - f^\ast \|_{\infty}$ by only assuming $f^\ast \in L^2(\m X)$. Reductions of the bound when $f^\ast$ is in either of the smoothness classes $\Theta_{\So}^{\alpha}$ or $\Theta_{\H}^{\alpha}$ are discussed subsequently. 

\begin{theorem}[Sup-norm bounds for KRR estimator]\label{Thm:KRR_bound}
Assume the kernel $K$ satisfies assumptions  {\bf (B)} and {\bf (E)}. Define $h$ via the equation $\lambda = h^{2\alpha}$ ($\alpha>1/2$). 
Then, with probability at least $1-n^{-10}$ with respect to the randomness in $\{(X_i,\,w_i)\}_{i=1}^n$, the following error estimate holds,
\begin{align}\label{Eqn:First_Order_bound}
 \|\wht{f}_{n, \lambda} -  f^\ast\|_\infty \leq 2\|P_\lambda f^\ast\|_\infty  + C \, \sigma\, \sqrt{\frac{\log n}{nh}}.
\end{align} 
Moreover, with the same probability, the following higher-order expansion holds,
\begin{equation}\label{Eqn:Second_Order_bound}
\bigg\|\wht{f}_{n, \lambda} -  F_\lambda f^\ast -\frac{1}{n}\,\sum_{i=1}^n w_i\,\wt{K}_{X_i} \bigg\|_\infty 
\leq C' \, \gamma_n \,\bigg(2\|P_\lambda f^\ast\|_\infty  + C \, \sigma\, \sqrt{\frac{\log n}{nh}} \,\bigg),
\end{equation}
where 
\begin{equation}\label{eq:gamma_n}
\gamma_n = \max\Big\{ 1, \,n^{-1+1/(2\alpha)}\,h^{-1/(2\alpha)}\, \sqrt{\log n }\Big\} \, \sqrt{\frac{\log n}{nh}}. 
\end{equation}
Here, $C$ and $C'$ are constants independent of $(n,\,\lambda,\,h,\,\sigma)$, and $\wt{K}$ is the equivalent kernel of $K$ defined in \S \ref{subsec:EK}.
\end{theorem}

\begin{remark}[Sub-Gaussian errors]\label{rem:subg}
An inspection of the proof of Theorem \ref{Thm:KRR_bound} reveals that we haven't made explicit use of the normality of the $w_i$s. Indeed, the conclusions of the theorem, and all subsequent results, extend to sub-Gaussian errors. 
\end{remark}

The first-order bound in \eqref{Eqn:First_Order_bound} has two components: a bias term, $\|P_{\lambda} f^\ast\|_{\infty}$, and a variance term, $\sigma \, (\log n)^{1/2} \, (nh)^{-1/2}$. Since $P_{\lambda} f^\ast = f^\ast - F_{\lambda} f^\ast$, $\|P_{\lambda} f^\ast\|_{\infty}$ measures the closeness (in terms of the supremum norm) between the true function $f^\ast$ and its convolution with the equivalent kernel $\wt{K}$, $F_{\lambda} f^\ast$. Recall also that $F_{\lambda} f^\ast$ is the solution to the population-level score equation. The smoothing parameter $\lambda$ provides the trade-off between the bias and variance; larger $\lambda$ (stronger regularization) reduces variance at the cost of increasing the bias, and vice versa. Notice that a typical analysis of the KRR estimator using basic inequality (such as \cite{yang2017randomized}) requires a lower bound on the regularization parameter $\lambda$, while our analysis is free of such an assumption.
Under additional assumptions on $f^\ast$, more explicit bounds can be obtained for the bias term $\|P_{\lambda} f^\ast\|_{\infty}$. For example, if $f^\ast \in \Theta_{\So}^{\alpha}(B)$, then one can show that $\|P_{\lambda} f^\ast\|_{\infty} \lesssim h^{\alpha - 1/2}$. Similarly, for $f^\ast \in \Theta_{\H}^{\alpha}(B)$, one obtains the bound $\|P_{\lambda} f^\ast\|_{\infty} \lesssim h^{\alpha}$. Choosing $h$ optimally in either situation lead to the following explicit bounds. 
\begin{corollary}[Minimax rates for Sobolev and H\"{o}lder classes]\label{cor:krr_supnorm}
Suppose $f^\ast \in \Theta_{\So}^{\alpha}(B)$ defined in \eqref{eq:Sob}. Set 
$$
h = \bigg( \frac{ B^2 n}{ \sigma^2 \log n} \bigg)^{-1/(2 \alpha)}. 
$$
Then, with probability at least $1-n^{-10}$ with respect to the randomness in $\{(X_i,\,w_i)\}_{i=1}^n$,
\begin{align}
 \|\wht{f}_{n, \lambda} -  f^\ast\|_\infty \lesssim B^{1/(2 \alpha)} \, \bigg( \frac{\sigma^2 \log n}{n} \bigg)^{\frac{(\alpha - 1/2)}{2 \alpha}}. 
\end{align}
Next, suppose $f^\ast \in \Theta_{\H}^{\alpha}(B)$ defined in \eqref{eq:Hol}. Set 
$$
h = \bigg( \frac{ B^2 n}{ \sigma^2 \log n} \bigg)^{-1/(2 \alpha + 1)}. 
$$
Then, with probability at least $1-n^{-10}$ with respect to the randomness in $\{(X_i,\,w_i)\}_{i=1}^n$,
\begin{align}
 \|\wht{f}_{n, \lambda} -  f^\ast\|_\infty \lesssim B^{1/(2 \alpha + 1)} \, \bigg( \frac{\sigma^2 \log n}{n} \bigg)^{\frac{\alpha}{(2 \alpha + 1)}}. 
\end{align}
\end{corollary}
Corollary \ref{cor:krr_supnorm} implies that the rate of convergence (up to logarithmic terms) of the KRR estimator to the truth in supremum norm is $n^{-(\alpha - 1/2)/(2 \alpha)}$ and $n^{-\alpha/(2 \alpha + 1)}$ for $\alpha$-smooth Sobolev and H\"{o}lder functions respectively. For functions in the Sobolev class, the KRR estimator concentrates at the usual minimax rate of $n^{-\alpha/(2 \alpha + 1)}$ under the $L_2$ norm \cite{caponnetto2007optimal,zhang2005learning}. However, it is known \cite{butucea2001exact,brown1996constrained} that under the point-wise and/or supremum norm, the minimax rate for $\alpha$-smooth Sobolev functions deteriorates to $n^{-(\alpha - 1/2)/(2 \alpha)}$. Hence, the KRR estimator achieves the minimax rate under the supremum norm as well. 
For the H\"{o}lder class, the minimax rate remains the same under the $L_2$ and $L_{\infty}$ norms, and the KRR estimator achieves the minimax rate. 

The higher-order expansion in the display \eqref{Eqn:Second_Order_bound} provides a finer insight into the distributional behavior of $\wht{f}_{n, \lambda}$. Let $U(\cdot)$ denote the zero-mean random process $\sqrt{h/n} \sum_{i=1}^n w_i \wt{K}_{X_i}$ with covariance function $\wht{C}_n$, so that $\wht{C}_n(x, x') = \bbE U(x) U(x')$. For any fixed $x$, the law of $U(x)$ can be approximated by the law of a $\mbox{N}(0, \wht{C}_n(x))$ distribution by the central limit theorem, and \eqref{Eqn:Second_Order_bound} can be used to establish that the law of $\sqrt{nh} (\wht{f}_{n, \lambda}(x) - F_{\lambda} f^\ast(x))$ is close to a $\mbox{N}(0, \wht{C}_n(x))$ distribution. Indeed, we shall establish a stronger result that the law of the process $\sqrt{nh} (\wht{f}_{n, \lambda}(\cdot) - F_{\lambda} f^\ast(\cdot))$ can be approximated by a Gaussian process $\mbox{GP}(0, \wht{C}_n)$ in the Kolmogorov distance. The point-wise and uniform approximation results will be crucial to prove asymptotic validity of Bayesian point-wise and simultaneous credible bands in \S \ref{sec:Bayes_theory}. 

Theorem \ref{Thm:KRR_bound} specializes to the noiseless KRR problem in a straightforward fashion upon setting $\sigma = 0$ in the bounds \eqref{Eqn:First_Order_bound} and \eqref{Eqn:Second_Order_bound}. As noted in \S \ref{sec:gp_reg}, the posterior variance function can be expressed as the solution to a noiseless KRR problem, motivating our interest in such situations. The following corollary states a general result for the noiseless case which is used subsequently to analyze the posterior variance function. 

\begin{corollary}[Sup-norm bounds for KRR estimator: noiseless case]\label{cor:KRR_bound}
Consider a noiseless version of the KRR problem in \eqref{eq:krr} where $Y_i = f^\ast(X_i)$. Assume the kernel $K$ satisfies assumptions  {\bf (B)} and {\bf (E)}. Define $h$ via the equation $\lambda = h^{2\alpha}$. 
Then, with probability at least $1-n^{-10}$ with respect to the randomness in $\{X_i\}_{i=1}^n$, the following error estimate holds,
\begin{align}\label{Eqn:fob1}
 \|\wht{f}_{n, \lambda} -  f^\ast\|_\infty \leq 2\|P_\lambda f^\ast\|_\infty.
\end{align} 
Moreover, with the same probability, 
\begin{equation}\label{Eqn:sob1}
\|\wht{f}_{n, \lambda} -  F_\lambda f^\ast \|_\infty = \|(f^\ast - \wht{f}_{n, \lambda}) - P_\lambda f^\ast \|_\infty 
\leq C' \, \gamma_n \|P_\lambda f^\ast\|_\infty,
\end{equation}
where $\gamma_n$ is as in \eqref{eq:gamma_n}. As before, $C$ and $C'$ are constants independent of $(n,\,\lambda,\,h)$, and $\wt{K}$ is the equivalent kernel of $K$ defined in \S \ref{subsec:EK}.
\end{corollary}

The proof of Corollary \ref{cor:KRR_bound} follows directly from tracking the proof of Theorem \ref{Thm:KRR_bound} with $\sigma = 0$, and hence omitted. Implications of Corollary \ref{cor:KRR_bound} for the posterior variance function are discussed in the next section.

\section{Convergence limit of Bayesian posterior}\label{sec:Bayes_theory}
We now use the sup-norm KRR bounds developed in \S \ref{sec:supnorm_krr} to analyze the GP posterior $f \mid \mb D_n \sim \mbox{GP}(\widehat{f}_n, \wt{C}_n^B)$ defined in \eqref{eq:meancov}. \\[1ex]
{\bf Implications for the posterior mean function:}  As noted in \S \ref{sec:gp_reg}, the posterior mean $\wht{f}_n$ under the $\mbox{GP}(0, \sigma^2 (n \lambda)^{-1} K)$ prior coincides with the KRR estimator $\wht{f}_{n, \lambda}$. Hence, the conclusions of Theorem \ref{Thm:KRR_bound} apply to $\wht{f}_n$, that is, with probability at least $1 - n^{-10}$, 
\begin{equation}\label{Eqn:mean_function_expansion}
\begin{aligned}
& \|\wht{f}_{n} -  f^\ast\|_\infty \leq 2\|P_\lambda f^\ast\|_\infty  + C \, \sigma\, \sqrt{\frac{\log n}{nh}},\qquad\mbox{and}\\
&\bigg\|\widehat{f}_n-  f^\ast -\bigg(\frac{1}{n}\,\sum_{i=1}^n w_i\,\wt{K}_{X_i} - P_\lambda f^\ast\bigg)\bigg\|_\infty  \leq C \, \gamma_n\,\bigg(\|P_\lambda f^\ast\|_\infty  +  \sigma\, \sqrt{\frac{\log n}{nh}} \,\bigg):\,=\delta_n,
\end{aligned}
\end{equation}

 In particular, for optimal choices of the prior precision parameter $\lambda$ as in Corollary \ref{cor:krr_supnorm}, the posterior mean achieves the minimax rates for the Sobolev and H\"{o}lder classes under the supremum norm. \\[1ex]
{\bf Implications for the posterior variance function:} Recall from \S~\ref{sec:gp_reg} that the posterior covariance function $\wt{C}^B_n$ admits the representation $\sigma^{-2}\,n \lambda\, \wt{C}^B_n(x,x')  =  K_{x}(x') - \wht{K}_{x}(x')$, where $\wht{K}_x$ is the solution to the noiseless KRR problem in \eqref{eq:var_krr}.   From  \eqref{Eqn:fob1} in  Corollary~\ref{cor:KRR_bound}, it follows with at least $1-n^{-10}$, 
\begin{align*}
\Big\|\sigma^{-2}\,n \lambda\, \wt{C}^B_n(x,\cdot) \Big\|_\infty \leq C \|P_\lambda \,K_x\|_\infty.
\end{align*}
A simple calculation yields 
$$
P_{\lambda} K_x(\cdot) = \sum_{j=1}^{\infty} (1 - \nu_j) \mu_j \, \psi_j(x) \psi_j(\cdot) = \lambda \sum_{j=1}^{\infty} \nu_j  \,\psi_j(x) \psi_j(\cdot) = \lambda \wt{K}_x(\cdot). 
$$
Since $\{\psi_j\}$ are uniformly bounded, we have 
\begin{align}\label{eq:nubound}
\sum_{j=1}^{\infty} \nu_j  \,\psi_j(x) \psi_j(x') \precsim \lambda^{-1/(2\alpha)}= h^{-1}.
\end{align}
From \eqref{Eqn:sob1} in Corollary~\ref{cor:KRR_bound}, we obtain that with probability at least $1-n^{-10}$, 
\begin{align*}
\Big\|\sigma^{-2}\,n \lambda\, \wt{C}^B_n(x,\cdot) - P_\lambda \,K_x\Big\|_\infty \leq C\,\gamma_n\,\|P_\lambda \,K_x\|_\infty.
\end{align*}
Combining with the previous display, 
\begin{align}\label{Eqn:Cov_diff}
\Big\|\sigma^{-2}\,n \, \wt{C}^B_n(x,\cdot) - \wt{K}_x\Big\|_\infty \leq C\,\gamma_n\,\|\wt{K}_x\|_\infty.
\end{align}
In particular, this relationship between the rescaled posterior covariance function $\sigma^{-2}\,n\,\wt{C}^B_n$ and the equivalent kernel function $\wt{K}$ leads to a practically useful way to numerically approximate the equivalent kernel function without explicitly conducting the eigen-decomposition to the original kernel function $K$.

Let $\wt{f}^B\sim \mbox{GP}(0,\wt{C}^B_n)$. Then the conditional distribution of $\wht{f}_{n} +\wt{f}^B$ given $\mb D_n$ is the posterior distribution of the mean function $f \mid \mb D_n$, or in other words, the law of $\wt{f}^B$ given 
$\mb D_n$ is that of the centered posterior distribution. When studying second-order properties of the posterior, it will be useful to work with a scaled version of the posterior with a $\sqrt{nh}$ scaling, which has the same covariance function as $\sqrt{nh} \, \wt{f}^B$; $\sqrt{nh} \, \wt{f}^B \mid \mb D_n \sim \mbox{GP}(0, nh \, \wt{C}^B_n)$. Write 
$$
nh \, \wt{C}^B_n = \sigma^2 h \, \big( \sigma^{-2} n \, \wt{C}^B_n \big). 
$$
The approximation error bound in \eqref{Eqn:Cov_diff} motivates us to define a related GP as $W^B\sim \mbox{GP}(0,\widehat{C}^B_n)$, with 
\begin{align*}
\widehat{C}^B_n(x,x'):\,=\sigma^2\, h \, \wt{K}(x, x') = \sigma^2\, h\,\sum_{j=1}^\infty \nu_j \psi_j(x)\,\psi_j(x'). 
\end{align*}
Notice that inequality~\eqref{Eqn:Cov_diff} implies that following sup-norm difference between covariance functions of $\sqrt{nh} \wt{f}^B$ and $W^B$,
\begin{align}\label{Eqn:Cov_approx_bound}
\sup_{x,x'}\Big| nh\, \wt{C}^B_n(x,\,x')  - \widehat{C}^B_n(x,\,x')\Big| \leq C\, \sigma^2 h \, \gamma_n\, \sup_{x,x'} \big|\wt{K}(x,\,x')\big| \lesssim \gamma_n,
\end{align}
since $\sup_{x, x'} | \wt{K}(x, x') | \lesssim \sum_{j=1}^{\infty} \nu_j \lesssim h^{-1}$ from  \eqref{eq:nubound}. Importantly, $\widehat{C}^B_n$ is a fixed function (depending only on the eigenbasis and the scaling $h$) unlike $\widetilde{C}_n^B$ which involves the random design $\{X_i\}_{i=1}^n$. We shall make repeated use of this approximation-error bound  in \eqref{Eqn:Cov_approx_bound} in the sequel for studying sup-norm posterior convergence rate and frequentist converge of Bayesian credible intervals/bands.

\subsection{Point-wise and sup norm posterior convergence rate}
We first state a result for the posterior rate of contraction around $f^\ast$ in terms of the point-wise and supremum norm. The Hellinger and total variation norms are by far the most common metrics for establishing posterior contraction rates and there is only a recent small literature for stronger norms \cite{gine2011rates,castillo2014bayesian,hoffmann2015adaptive,yoo2016supremum}. The following theorem establishes that the GP posterior concentrates around the true function in both the point-wise and supremum norms at the respective minimax rates for the Sobolev and H\"{o}lder classes. 

\begin{theorem} \label{thm:ratepw&sn}
Assume the kernel $K$ satisfies assumptions  {\bf (B)} and {\bf (E)}. Define $h$ via the equation $\lambda = h^{2\alpha}$. 
If $h\asymp \big\{B^2 n/(\sigma^2 \, \log n)\big\}^{-1/(2\alpha+1)}$ for $f^\ast \in\Theta_{\H}^{\alpha}(B)$ and $h \asymp \big\{B^2 n/(\sigma^2 \, \log n)\big\}^{-1/(2\alpha)}$ for $f^\ast \in\Theta_{\So}^{\alpha}(B)$, then, with probability at least $1-n^{-10}$ with respect to the randomness in $\{X_i, w_i\}_{i=1}^n$, the following error estimates for the squared point-wise risk hold for any $x^*$ in $\mathcal{X}$,
\begin{align}\label{Eqn:fob1}
 \bbE[ |f(x^\ast) - f^\ast(x^\ast)|^2 \mid \mathbb{D}_n] \precsim  
\begin{cases}
B^{2/(2\alpha+1)} \bigg( \frac{\sigma^2 \log n}{n} \bigg)^{\frac{2\alpha}{(2 \alpha + 1)}}  
& \text{if } f^\ast \in \Theta_{\H}^{\alpha}(B),  \\
B^{1/\alpha} \bigg( \frac{\sigma^2 \log n}{n} \bigg)^{\frac{2\alpha - 1}{2 \alpha }} 
& \text{if } f^\ast \in \Theta_{\So}^{\alpha}(B).
\end{cases}%
\end{align}
With the same probability and the same choice of $h$, the following error estimates for the squared supremum norm hold,
\begin{align}\label{Eqn:fob2}  
\bbE[ \|f- f^*\|_{\infty}^2 \mid \mathbb{D}_n] \precsim  
\begin{cases}
B^{2/(2\alpha+1)} \bigg( \frac{\sigma^2 \log n}{n} \bigg)^{\frac{2\alpha}{(2 \alpha + 1)}} 
& \text{if } f^\ast \in \Theta_{\H}^{\alpha}(B),  \\
B^{1/\alpha} \bigg( \frac{\sigma^2 \log n}{n} \bigg)^{\frac{2\alpha - 1}{2 \alpha }}
& \text{if } f^\ast \in \Theta_{\So}^{\alpha}(B).   
\end{cases}
\end{align} 
\end{theorem}

Two remarks are in order. First, we can also apply similar techniques to show error estimates regarding the derivatives $f^{(k)}$ of $f$ with order $k$ up to $\lfloor \alpha\rfloor$ by identifying the posterior covariance function of $f^{(k)}$ (whose posterior is also a GP) with certain noiseless KRR estimate, which also applies to results in the following subsections. Due to the space constraint, we omit these results.
Second, since the sup-norm is stronger than the $L^2$ norm, our result shows that by introducing the scaling $(n\lambda)^{-1}$ in the covariance kernel of the prior GP, we can eliminate the mismatch \cite{van2008rates,vaart2011information} between the smoothness level (which is $\alpha+1/2$ for the univariate case and $\alpha+d/2$ for the $d$-variate case) of the RKHS of the $\alpha$-Mat\'{e}rn kernel and the best smoothness level of the truth (which is $\alpha$ without this scaling) that the posterior can adapt to.

\subsection{Frequentist coverage of point-wise posterior credible intervals}

In the following, we leverage on the (higher order) expansions for both the posterior mean \eqref{Eqn:mean_function_expansion} and the variance \eqref{Eqn:Cov_approx_bound} to derive the frequentist coverage properties of point-wise Bayesian credible intervals. 

Let $\Phi$ denote the $\mbox{N}(0, 1)$ c.d.f. and for any $\gamma \in (0, 1)$, let $z_{\gamma}$ denote the $\gamma$th standard normal quantile with $\Phi(z_{\gamma}) = \gamma$. Since the posterior distribution of $f$ is GP$(\wht{f}_{n}, \,\wt{C}^B_n)$, we consider for any $x\in\m X$ a point-wise credible interval centered at $\wht{f}_{n}(x)$ with level $\beta\in(0,1)$ as 
$$\mbox{CI}_n(x;\,\beta) = \big[\wht{f}_{n}(x) - l_n(x;\,\beta),\, \wht{f}_{n}(x) + l_n(x;\,\beta)\big],$$
where the half length 
$$l_n(x;\,\beta) = z_{(1+\beta)/2} \, \sqrt{ \wt{C}^B_n(x,\, x) } \, ,$$ 
is chosen so that the posterior probability of $f(x)$ falling into the credible interval is $\beta$, or
\begin{align*}
\mb P\big[ f(x) \in \mbox{CI}_n(x;\,\beta) \,\big| \, \mb D_n\big] = \beta.
\end{align*}
Note that a combination of \eqref{Eqn:Cov_approx_bound} and the fact that $\sup_{x, x'} | \wt{K}(x, x') | \lesssim \sum_{j=1}^{\infty} \nu_j \lesssim h^{-1}$ implies that the size of the credible interval $\mbox{CI}_n(x;\,\beta)$ is $\m O((nh)^{-1/2})$, which is minimax-optimal under choices of $h$ in Theorem~\ref{thm:ratepw&sn}.
Our goal below is to investigate the frequentist coverage of $\mbox{CI}_n(x;\,\beta)$, and in particular, identify situations when 
\begin{align*}
\mb P_{\rho}\big[ f^\ast(x) \in \mbox{CI}_n(x;\,\beta) \,\big] \ge \beta.
\end{align*}
Note that $\mbox{CI}_n(x;\,\beta)$ is a random interval in the above display and $\mb P_{\rho}$ denotes the probability under the true data generating distribution that $\mbox{CI}_n(x;\,\beta)$ contains the true function value $f^\ast(x)$. Let $\gamma = (1 + \beta)/2$. Write 
\begin{align}
& \mb P_{\rho}\big[ f^\ast(x) \in \mbox{CI}_n(x;\,\beta) \,\big]  \notag \\
& = \mb P_{\rho} \bigg[ - z_{\gamma} \sqrt{nh \, \wt{C}^B_n(x,\, x) } \le \sqrt{nh} \big[\wht{f}_n(x) - f^\ast(x) \big] \le z_{\gamma} \sqrt{nh \, \wt{C}^B_n(x,\, x) } \, \bigg] \notag \\
& \stackrel{(i)} \approx \mb P_{\rho} \bigg[-z_{\gamma} \sqrt{\widehat{C}^B_n(x,\,x)} \le \sqrt{nh} \big[\wht{f}_n(x) - f^\ast(x) \big] \le z_{\gamma} \sqrt{\widehat{C}^B_n(x,\,x)} \, \bigg] \notag \\
& = \mb P_{\rho} \bigg[ \sqrt{nh} \big[\wht{f}_n(x) - F_{\lambda} f^\ast(x) \big] \in \Big[ \sqrt{nh} \, P_{\lambda} f^\ast(x) \pm z_{\gamma} \sqrt{\widehat{C}^B_n(x,\,x)}  \, \Big]  \bigg]. \label{eq:cov_illus1}
\end{align}

The approximation (i) in the above display follows from \eqref{Eqn:Cov_approx_bound} and the approximation bound will be made concrete inside the proof. Recall the process $U(\cdot) = \sqrt{h/n} \,\sum_{i=1}^n w_i \wt{K}_{X_i}(\cdot)$ from the discussion after Corollary \ref{cor:krr_supnorm}. Consider a GP $\widehat{W}\sim\mbox{GP}(0,\, \widehat{C}_n)$ whose covariance function $\widehat{C}_n$ matches the covariance 
function of the process $U(\cdot)$ over $\m X$, that is, for any $(x,\,x')\in\m X^2$,
\begin{align*}
\widehat{C}_n(x,x')=\mb E[U(x)U(x')]=\sigma^2\, h\,\mb E\big[\wt{K}(X, x)\,\wt{K}(X, x')\big] =\sigma^2\,  h\,\sum_{j=1}^\infty \frac{\psi_j(x)\,\psi_j(x')}{(1+h^{2\alpha}/\mu_j)^2}.
\end{align*}
The second display in~\eqref{Eqn:mean_function_expansion} can be used to establish that the law of the random process $\sqrt{nh}\,(\wht{f}_{n} - F_\lambda f^\ast)$ is (point-wise) close to that of $\widehat{W}$.  Substituting in \eqref{eq:cov_illus1} and tracking the approximation errors leads to the following theorem. 

\begin{theorem}[Frequentist coverage of posterior credible intervals]\label{thm:point-wise_cov}
There exists some constant $C$ independent of $(n,h)$ such that the frequentist coverage of $\mbox{CI}_n(x;\,\beta)$ satisfies that for any $x\in \m X$,
\begin{align*}
\Big| &\mb P\big[f^\ast(x) \in  \mbox{CI}_n(x;\,\beta) \big] - \big[ \Phi \big(u_n(x;\,\beta) + b_n(x)\big)  - \Phi \big(-u_n(x;\,\beta) + b_n(x)\big)\big]   \Big| \leq C\, \Big( \frac{1}{\sqrt{nh}} + \gamma_n+\delta_n\Big),
\end{align*}
where $u_n(x;\,\beta) = \sqrt{\widehat{C}^B_n(x,\,x)/\widehat{C}_n(x,x)}\, z_{(1+\beta)/2}$ is the inflated quantile and \\ $b_n(x) = \big\{\widehat{C}_n(x,x)\big\}^{-1/2}\,\sqrt{nh}\, P_\lambda f^\ast(x)$ is the bias at $x$.
\end{theorem}

We briefly comment on the source of each of the three approximation error terms on the right hand side of the previous display; details can be found inside the proof.  The distributional approximation of $\big\{\widehat{C}_n(x,x)\big\}^{-1/2} U(x)$ with a $\mbox{N}(0, 1)$ random variable incurs an error of $(nh)^{-1/2}$ from the Berry--Essen theorem.  The approximation of the appropriately standardized posterior mean $\big\{\widehat{C}_n(x,x)\big\}^{-1/2}\sqrt{nh}(\hat{f}_n(x)  -  F_\lambda f^*(x))$ with $\big\{\widehat{C}_n(x,x)\big\}^{-1/2} U(x)$ results in an error of $\delta_n$ from \eqref{Eqn:mean_function_expansion}.   The approximation  \eqref{Eqn:Cov_approx_bound} contributes an error of $\gamma_n$.  

Come back to the concrete case when $[0,1]$ is the unit interval and $K$ is the Mat\'{e}rn kernel with smoothness index $\alpha>1/2$, where $\psi_{2j}(x) = \cos(\pi j x)$, $\psi_{2j-1}=\sin(\pi jx)$, and $\mu_{2j}= \mu_{2j+1}\asymp j^{-2\alpha}$ for $j=0,1,\ldots$. Applying the identity $\cos(x-y)=\cos(x)\cos(y)+\sin(x)\sin(y)$, we can simplify the two covariance functions $\widehat{C}^B_n$ and $\widehat{C}_n$ into
\begin{align}\label{Eqn:C_matern}
\widehat{C}^B_n(x,x') =\sigma^2\,  h\,\sum_{j=1}^\infty \frac{\cos(\pi j (x-x'))}{1+h^{2\alpha}/\mu_{2j}}\quad \mbox{and}\quad \widehat{C}_n(x,x') =\sigma^2\,  h\,\sum_{j=1}^\infty \frac{\cos(\pi j (x-x'))}{(1+h^{2\alpha}/\mu_{2j})^2}.
\end{align}

Define\begin{align*}
C_{IR} :=&\, 
\widehat{C}^B_n(x,\,x)/\widehat{C}_n(x,x) = \sum_{j=1}^\infty \frac{1}{1+h^{2\alpha}/\mu_{2j}}\, \big/\, \sum_{j=1}^\infty \frac{1}{(1+h^{2\alpha}/\mu_{2j})^2},
\end{align*}
where we have applied identities in \eqref{Eqn:C_matern} so that the $C_{IR}$ is independent of $x$.  Clearly $C_{IR} > 1$.  
 For simplicity, we consider  
$f^\ast \in \Theta_{\H}^{\alpha_0}(B)$.  All  the results naturally extends to $f^\ast \in \Theta_{\So}^{\alpha_0}(B)$.
\begin{corollary}[Mat\'{e}rn kernel class] \label{cor:matern}
{\bf (Under-smooth)}  Suppose $f^\ast \in \Theta_{\H}^{\alpha_0}(B)$ for $\alpha_0 > \alpha$, then for any $x\in\m X$, $b_n(x) \to 0$ as $n\to\infty$, and
\begin{align*}
\mb P\big[f^\ast(x) \in  \mbox{CI}_n(x;\,\beta) \big] \to 2\,\Phi(C_{IR}\, z_{(1+\beta)/2}) -1 \in(\beta, 1), \quad\mbox{as $n\to\infty$}.
\end{align*}
{\bf (Smooth-match)} For any $\widetilde{\beta}\in(0,1)$, there is some sufficiently large constant $B>0$ and a sequence of functions $\{f^\ast_n\}$, each belongs to $\Theta_{\H}^{\alpha_0}(B)$, such that
\begin{align*}
\mb P\big[f_n^\ast(x) \in  \mbox{CI}_n(x;\,\beta) \big] \to \widetilde{\beta}, \quad\mbox{as $n\to\infty$}.
\end{align*}
{\bf (Over-smooth)}
 For any $\alpha_0<\alpha$, there always exists a function $f^\ast_{bad} \in \Theta_{\So}^{\alpha_0}(B)$, such that for any $x\in\m X$, $b_n(x) \to \infty$ as $n\to\infty$, and 
\begin{align*}
\mb P\big[f^\ast_{bad}(x) \in  \mbox{CI}_n(x;\,\beta) \big] \to 0, \quad\mbox{as $n\to\infty$}.
\end{align*}
\end{corollary}

\subsection{Frequentist coverage of simultaneous posterior credible bands}

In this subsection, we study the frequentist coverage of the following posterior credible band centered at the posterior mean $\wht{f}_{n}$ with level $\beta\in(0,1)$,
\begin{align*}
\mbox{CB}_n(\beta) = \Big\{f\in L^2(\m X):\, \big\|f - \wht{f}_n\big\|_\infty \leq r_n(\beta)\Big\},
\end{align*}
where the half length $r_n(\beta)$ is chosen so that posterior probability of $f$ falling into the credible band is $\beta$, i.e.,
\begin{align*}
\mb P\big[ f \in \mbox{CB}_n(\beta)\, \big|\, \mb D_n \big]= \beta.
\end{align*}
Whereas the point-wise intervals $\mbox{CI}_n(x; \beta)$ permitted an explicit description using the Gaussianity of $f(x) \mid \mb D_n$ for any $x$, the band $\mbox{CB}_n(\beta)$ needs to be defined implicitly due to the lack of a similar exact distributional result. Establishing the frequentist validity of the band $\mbox{CB}_n(\beta)$ follows a conceptually similar route as before; (i) approximate the sampling distribution of the standardized posterior mean $\sqrt{nh} \, (\wht{f}_n - F_{\lambda} f^\ast)$ by the GP $\wht{W} \sim \mbox{GP}(0, \wht{C}_n)$ define previously, and (ii) approximate the centered and scaled posterior measure $\sqrt{nh} \, (f - \wht{f}_n) \mid \mb D_n$ by the GP $\wht{W}^B \sim \mbox{GP}(0, \wht{C}_n^B)$. However, substantial work is necessary to obtain uniform counterparts of the point-wise approximations obtained previously. 

We first discuss approximation of the sampling distribution of $\sqrt{nh} \, (\wht{f}_n - F_{\lambda} f^\ast)$ in supremum norm. Recall that $\widehat{W}$ is defined as a GP with law $\mbox{GP}(0,\, \widehat{C}_n)$, where $\widehat{C}_n$ is the covariance function of the random process $U(\cdot)=\sqrt{h/n}\,\sum_{i=1}^n w_i\,\wt{K}_{X_i}(\cdot)$ over $\m X$, which is the leading term in the expansion of $\sqrt{n\,h}\,\big(\widehat{f}_n-F_\lambda\, f^\ast\big)$ in the second display in~\eqref{Eqn:mean_function_expansion}. Let $\widehat{Z}_n^M=\|\widehat{W}\|_\infty$ and $Z_n^M=\big\|U\big\|_\infty$. We will first show by applying the results in \cite{chernozhukov2014gaussian} on Gaussian approximation to the suprema of empirical processes that the distributions of $\widehat{Z}_n^M$ and $Z_n^M$ are close with respect to the Kolmogorov distance. 

We also define $\widetilde{Z}_n^M=\sqrt{n\,h}\,\big\|\widehat{f}_n-F_\lambda\, f^\ast\big\|_\infty$ as the rescaled sup-norm deviation of the posterior mean function from its population counterpart $F_\lambda f^\ast = (I-P_\lambda)f^\ast$. By applying an anti-concentration bound \citep{chernozhukov2014anti} for GP, we obtain by combining the second display in~\eqref{Eqn:mean_function_expansion} that the distribution of $\widetilde{Z}_n^M$ can be well-approximated by that of $Z_n^M$. Finally, a combining of these two approximation results implies that the distribution of $\widetilde{Z}_n^M$ can be well-approximated by
$\widehat{Z}_n^M$, the supremum of a Gaussian process $\widehat{W}$, as summarized in the following theorem.

\begin{theorem}[Gaussian process approximation of the posterior mean function]\label{Thm:GP_mean_function}
There exists some constant $C$ independent of $(n,h)$ such that for any $t\geq 0$,
\begin{align}\label{Eqn:PosteriorModeGP}
\Big|\mb P\big[\widehat{Z}_n^M \leq t\big] - \mb P\big[\wt{Z}_n^M \leq t\big] \Big| \leq C\, \frac{\sqrt{\log n}}{(nh)^{1/8}}.
\end{align}
\end{theorem}

Next, we show that the posterior measure can also be uniformly approximated by another specially designed GP. Recall that $\sqrt{nh} \, \wt{f}^B \sim \mbox{GP}(0, nh \, \widetilde{C}^B_n)$, where $\wt{f}^B$ has the law of $(f - \wht{f}_n) \mid \mb D_n$. The covariance function $\wt{C}^B_n$ depends on the random design $\{X_i\}_{i=1}^n$ and is hard to directly work with. For this reason, we defined a population level GP
$W^B\sim \mbox{GP}(0,\widehat{C}^B_n)$ whose covariance function $\widehat{C}^B_n$ provides good approximation \eqref{Eqn:Cov_approx_bound} to $nh \, \wt{C}^B_n$ in sup-norm. Let $\widetilde{Z}^B_n=\sqrt{n\,h}\,\|\wt{f}^B\|_\infty$ and $Z^B_n =\|W^B\|_\infty$. In the next theorem, we show that the distributions of $\widetilde{Z}^B_n$ and its population level counterpart $Z^B_n$ are close with respect to the Kolmogorov distance. 
To prove Theorem \ref{Thm:GP_post}, we develop a new Gaussian comparison inequality; see Theorem \ref{Thm:GP_comparison} in \S \ref{Sec:Proof_Thm_GP_post}; that explicitly bounds the Kolmogorov distance between two GPs in terms of the supremum norm difference between their covariance functions. The comparison inequality extends the Sudakov--Fernique inequality for finite-dimensional Gaussians as stated in \cite{chatterjee2005error} to GPs.

\begin{theorem}[Gaussian process approximation of the centered posterior measure]\label{Thm:GP_post}
There exists some constant $C$ independent of $(n,h)$ such that for any $t\geq 0$,
\begin{align}\label{Eqn:PosteriorLimit_GP}
\Big|\mb P\big[\widetilde{Z}^B_n \leq t \, \big| \, \mb D_n\big] - \mb P\big[Z^B_n \leq t\big] \Big| \leq C\, \gamma_n^{1/3} \,\log n.
\end{align}
\end{theorem}

From Theorem~\ref{Thm:GP_mean_function} and Theorem~\ref{Thm:GP_post}, we only need to compare $\widehat{W}\sim\mbox{GP}(0,\, \widehat{C}_n)$ and $W^B\sim \mbox{GP}(0,\widehat{C}^B_n)$ in order to study the frequentist coverage of the Bayesian credible band $\mbox{CB}_n(\beta)$. Let $q_n^B(\beta)$ denote the $\beta$-th quantile of the population level random variable $Z^B_n =\|W^B\|_\infty$.

\begin{theorem}[Frequentist coverage of posterior credible bands]\label{thm:band_cov}
There exists some constant $C$ independent of $(n,h)$, such that for any $\beta\in(0,1)$,
\begin{align}\label{Eqn:quantile_approx}
\Big| r_n(\beta) - \frac{1}{\sqrt{nh}}\, q_n^B(\beta)\Big| \leq \frac{C}{\sqrt{nh}}\, \gamma_n^{1/3}\,\log n.
\end{align}
Moreover, if the bias term $\|P_\lambda f^\ast\|_\infty$ satisfies $\sqrt{nh}\,\|P_\lambda f^\ast\|_\infty\to 0$ as $n\to\infty$, then 
\begin{align}\label{Eqn:small_bias}
\Big|\mb P\big[ f^\ast \in \mbox{CB}_n(\beta)\big] - \mb P\big[ \|\widehat{W}\|_\infty \leq q_n^B(\beta)\big]\Big| \leq C\, \Big(\gamma_n^{1/3}\,\log n  + \frac{\sqrt{\log n}}{(nh)^{1/8}} + \sqrt{nh}\,\|P_\lambda f^\ast\|_\infty\Big) \to 0,
\end{align}
as $n\to\infty$ and $h\to 0$, where $\mb P\big[ \|\widehat{W}\|_\infty \leq q_n^B(\beta)\big] \in (\beta, 1)$.
On the other hand, if the bias term $\|P_\lambda f^\ast\|_\infty$ satisfies $\sqrt{nh}\,\|P_\lambda f^\ast\|_\infty\to \infty$ as $n\to\infty$, then 
\begin{align}\label{Eqn:large_bias}
\mb P\big[ f^\ast \in \mbox{CB}_n(\beta)\big] \to 0,\quad\mbox{as $n\to\infty$ and $h\to 0$.}
\end{align}
\end{theorem}

Since $(nh)^{-1/2}\, q_n^B(\beta)= \m O((nh)^{-1/2})$, the width of the simultaneous credible band is again minimax-optimal under choices of $h$ in Theorem~\ref{thm:ratepw&sn}.
Unlike the point-wise case where $\sqrt{nh}\,(\widehat{f}_n(x)-f^\ast(x))$ weakly converges to a nondegenerate distribution, in the simultaneous case, the two $n$-dependent population level GPs $W^B$ and $\widehat{W}$ generally do not weakly converge to non-degenerate laws (tight GPs) over $[0,1]$ as $n\to \infty$ and $h\to 0$, due to the kernel type estimator form of the process $U$ with an $(n,\,h)$ dependent kernel (see, for example, \cite{nishiyama2011impossibility,stupfler2016weak}). 

Similar to Corollary~\ref{cor:matern}, we have the following corollary when specialized to the Mat\'{e}rn kernel with smoothness index $\alpha>1/2$ over $[0,1]$, and omit the proof.
\begin{corollary}[Mat\'{e}rn kernel class] \label{cor:matern_sim}
{\bf (Under-smooth)}  Suppose $f^\ast \in \Theta_{\H}^{\alpha_0}(B)$ for $\alpha_0 > \alpha$, then the bias term $\sqrt{nh}\,\|P_\lambda f^\ast\|_\infty\to 0$ as $n\to\infty$, and
\begin{align*}
\mb P\big[ f^\ast \in \mbox{CB}_n(\beta) \big] \in(\beta, 1), \quad\mbox{as $n\to\infty$}.
\end{align*}
{\bf (Smooth-match)} For any $\widetilde{\beta}\in(0,1)$, there is some sufficiently large constant $B>0$ and a sequence of functions $\{f^\ast_n\}$, each belongs to $\Theta_{\H}^{\alpha_0}(B)$, such that
\begin{align*}
\mb P\big[f^\ast \in \mbox{CB}_n(\beta)\big] \to \widetilde{\beta}, \quad\mbox{as $n\to\infty$}.
\end{align*}
{\bf (Over-smooth)}
 For any $\alpha_0<\alpha$, there always exists a function $f^\ast_{bad} \in \Theta_{\H}^{\alpha_0}$ such that the bias term $\sqrt{nh}\,\|P_\lambda f^\ast\|_\infty\to \infty$ as $n\to\infty$, and 
\begin{align*}
\mb P\big[ f^\ast \in \mbox{CB}_n(\beta)\big] \to 0, \quad\mbox{as $n\to\infty$}.
\end{align*}
\end{corollary}
%
%

\section{Simulation study}
In the following,  we numerically investigate the behavior of the point-wise and simultaneous credible intervals for a certain $f^*$ when  a Gaussian process prior with Mat\'{e}rn covariance kernel given by 
\begin{eqnarray}\label{eq:matern}
K(x, y) = \frac{2^{1-\nu}}{\Gamma(\nu)} \big\{\sqrt{2\nu} \, |x-y| \big\}^\nu B_\nu\big\{ \sqrt{2\nu} \,  |x-y|\big\}, \quad x,y \in [0, 1],
\end{eqnarray}
where $B_\nu$ is the modified Bessel function of the second kind for $0 < \nu < \infty$.  We recall that $K$ with $\nu=(\alpha - 1/2)$ satisfies {\bf (B)} and {\bf (E)} when expanded with respect to the Fourier basis 
$\psi_{2j-1}(x) = \sin(\pi j x), \ \psi_{2j}(x) = \cos(\pi j x), \ j = 1, 2, \ldots$. Also the eigen-values satisfy $\lambda_{2j-1}  =\lambda_{2j} \asymp j^{-2\alpha}$.  

\begin{figure}[htbp!]
\centering
\begin{subfigure}{.95\linewidth}
    \includegraphics[width=\linewidth]{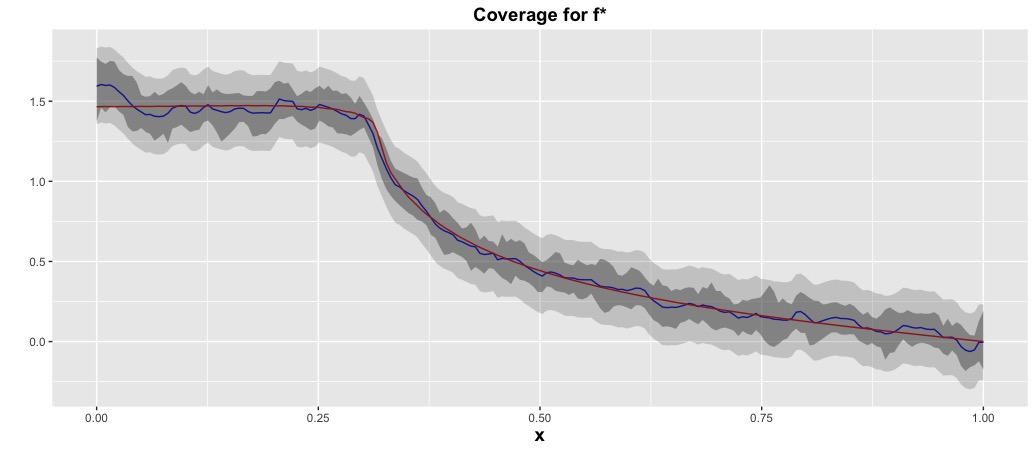}
    \caption{$n=200$}
    \label{fig11}
\end{subfigure} 
\begin{subfigure}{.95\linewidth}
    \includegraphics[width=\linewidth]{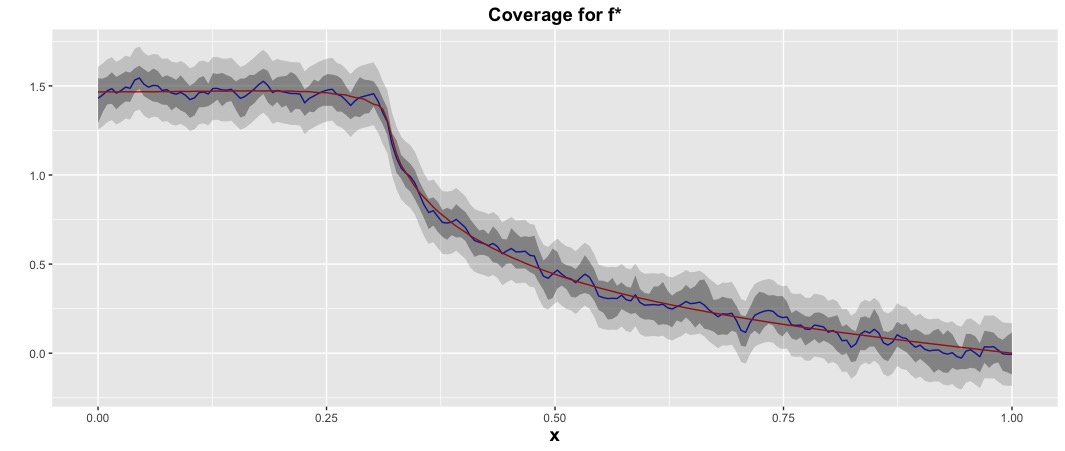}
    \caption{$n=500$}
        \label{fig12}
\end{subfigure}
\begin{subfigure}{.95\linewidth}
    \includegraphics[width=\linewidth]{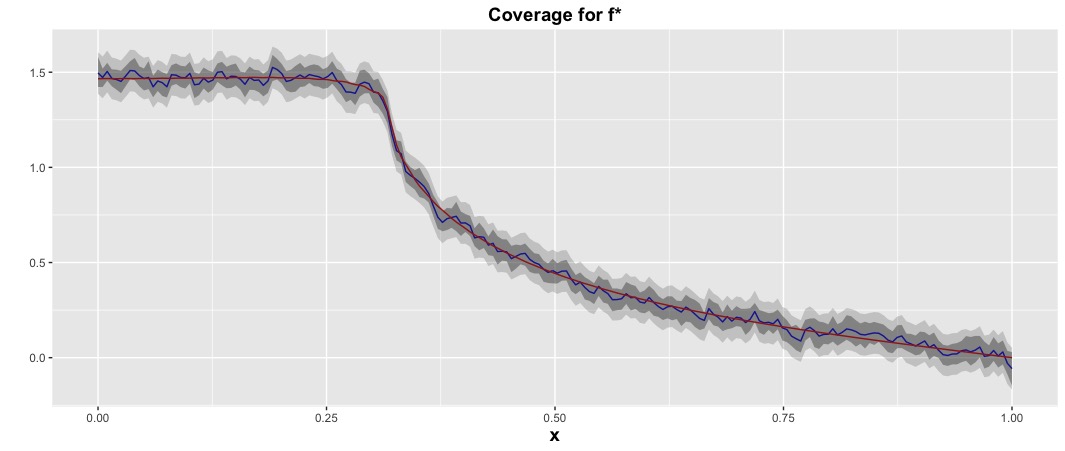}
    \caption{$n=2000$}
        \label{fig13}
\end{subfigure}
\caption{Credible intervals with Mat\'{e}rn with $\nu = 0.1$. The red and the blue lines represent $f^*$  and the point-wise posterior mean respectively.  The upper and lower brackets of the dark-grey and the light-grey shaded regions represent the point-wise and the simultaneous $95\%$ credible intervals.  
 }
 \label{fig:p}
\end{figure}

\begin{figure}[htbp!]
\centering
\begin{subfigure}{.95\linewidth}
    \includegraphics[width=\linewidth]{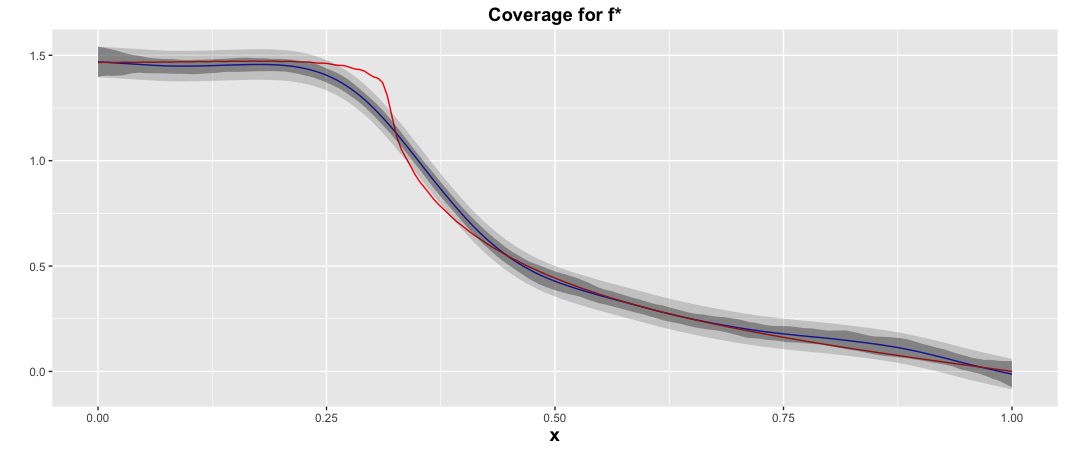}
    \caption{$n=200$}
    \label{fig21}
\end{subfigure} 
\begin{subfigure}{.95\linewidth}
    \includegraphics[width=\linewidth]{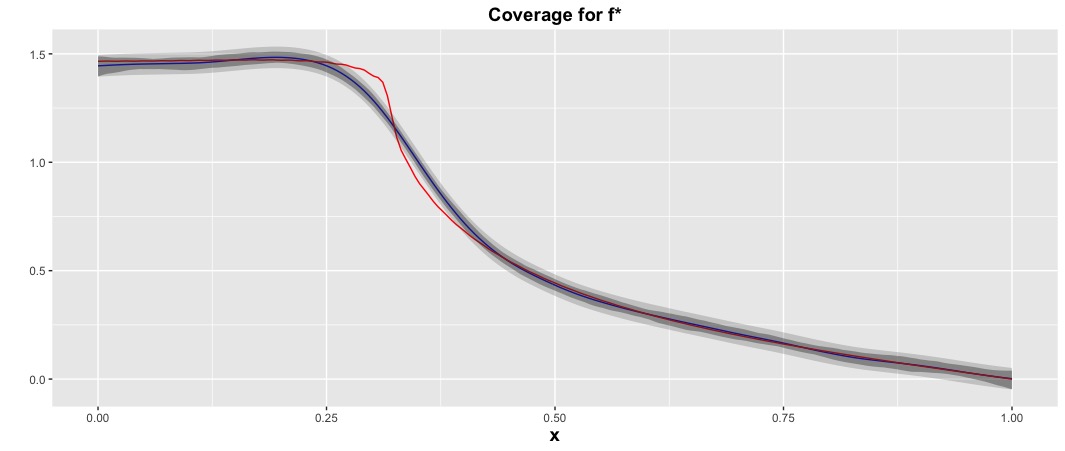}
    \caption{$n=500$}
        \label{fig22}
\end{subfigure}
\begin{subfigure}{.95\linewidth}
    \includegraphics[width=\linewidth]{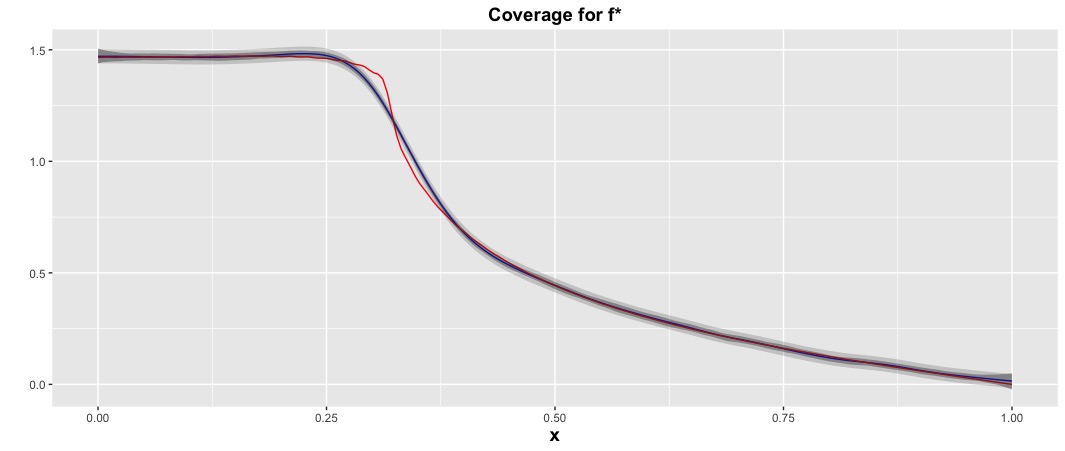}
    \caption{$n=2000$}
        \label{fig23}
\end{subfigure}
\caption{Credible intervals with Mat\'{e}rn with $\nu = 1.2$. The red and the blue lines represent $f^*$  and the point-wise posterior mean respectively.  The upper and lower brackets of the dark-grey and the light-grey shaded regions represent the point-wise and the simultaneous $95\%$ credible intervals.  
 }
 \label{fig:n}
\end{figure}

We let  $f^*(x) =  \sum_{j=1}^\infty \{j^{-1.7} \sin j\} \cos(\pi(j-0.5)x)$.   Training data of size $n$ are drawn from the model $y_i=f^*(x_i)+ \epsilon_i$ with $\epsilon_i \sim \mbox{N}(0, \, 0.1^2)$ and $x_i \sim \mbox{Unif}(0, 1)$.  $200$ equally spaced points on $[0, 1]$ are chosen as the test $x$ values.   For Figure \ref{fig:p}, the prior has smoothness $0.6$ where  $f^* \in \Theta_{\H}^{\alpha_0}(B)$ for $\alpha_0  <0.7$.   For Figure \ref{fig:n}, the prior has smoothness $1.7$  and $f^*$ satisfies the counter example in the proof of Corollary \ref{cor:matern} in the over-smooth case.  The posterior mean and covariance functions are directly obtained using \eqref{eq:postme} and \eqref{eq:postcov}, from  where the point-wise $95\%$ credible intervals at the $200$ test points are constructed.  To obtain the simultaneous credible bands, $1000$ random samples are drawn from the posterior distribution of  the function evaluated at the $200$ test points. Then $r_n(0.95)$ is estimated from the 95\%  quantile of the posterior samples of $\| f - \hat{f}_{n, \lambda}\|_\infty$.  In the over-smooth and the under-smooth case,  it is suggestive from Figures \ref{fig:p} and \ref{fig:n} that  the simultaneous and  point-wise coverage tends to a non-zero fraction and zero in the respective cases as the sample size  increases to $\infty$, which is consistent with the prediction from our theory.

The coverage of the simultaneous $\beta \%$-credible intervals is  further investigated through a replicated simulation study.  
We consider $\beta = 0.80$ and $0.90$.  Using the same $f^*$, we consider three smoothness parameters $(0.1, 0.15, 1.2)$ of the Mat\'{e}rn covariance kernel.  $1000$ datasets were drawn from the same model and the proportion of datasets for which $f^*$ lies in the simultaneous interval  is recorded in Table \ref{tab:13}.  In the over-smooth case, the coverage tends to $0$ and in the under-smooth case, the coverage tends to a number between $\beta$ and $1$.  It may appear from the simulation study that in the under-smooth case, the coverage probability is tending to 1 as the sample size increases. This is primarily because $\nu$ is close to zero and the actual coverage probability becomes very close to 1 in these cases. As a result when the sample size is very large,  one needs a huge number of replicates to actually sample from the tail event that a dataset for which the simultaneous credible band does not cover the true function.   
\begin{table}[htbp]
\centering
\scalebox{1.05}{
\begin{tabular}{ccccccc} \toprule
 {\bf n} & \multicolumn{2}{c}{{\bf 200}} & \multicolumn{2}{c}{{\bf 500}} &\multicolumn{2}{c}{{\bf 2000}}  \\
 \cmidrule(lr){2-3} \cmidrule(lr){4-5} \cmidrule(lr){6-7}  \\
 \text{Prior} & {\bf 0.80} & {\bf 0.90} &  {\bf 0.80} & {\bf 0.90} & {\bf 0.80} & {\bf 0.90} \\ \hline 
  \text{Mat\'{e}rn}($\nu = 0.1$) & 0.977  & 0.993 & 0.995  & 0.999  & 0.998 & 0.999   \\
      \text{Mat\'{e}rn}($\nu = 0.15$) & 0.875  & 0.939 & 0.926 & 0.953  & 0.978& 0.993   \\
        \text{Mat\'{e}rn}($\nu = 1.2$) & 0 & 0 & 0 & 0  & 0 & 0  \\ \hline 
 \end{tabular}}
\caption{Simultaneous coverage probability over 1000 replicates.}
\label{tab:13}
\end{table}

\newpage
\section{Proofs}\label{sec:proof}
In this section, we provide proofs of the results in the paper.
For notational brevity, we drop $\lambda$ from the subscript of $\|\cdot\|_{\lambda}$ within the proofs; all instances of $\|\cdot\|$ refer to the RKHS norm (for the equivalent kernel) $\|\cdot\|_{\lambda}$. 
\subsection{Summary of notations}

For the reader's convenience, Table \ref{tab:notes} provides a summary of some of the key notations introduced earlier that are heavily used inside the proofs. 

\begin{table}[h!]
\caption{Notation reference} \label{tab:notes}
\begin{center}
\begin{tabular}{cl|l}
& \textbf{Symbol} & \textbf{definition} \\
\hline
& $\mathbb{D}_n$         & data, $\{(Y_i, X_i), i = 1, \ldots, n\}$ \\[0.5ex]
& $f^\ast$ & true function \\[0.5 ex]
& $w_i$ & Gaussian error, $w_i = Y_i - f^\ast(X_i) \sim \mbox{N}(0, \sigma^2)$\\[0.5 ex]
& $\wht{f}_n$   & posterior mean function, $\mb E [f \mid \mathbb{D}_n]$\\[0.5ex]
& $\wt{C}_n^B$ & posterior covariance function, $\wt{C}_n^B(x, x') = \mbox{Cov}(f(x), f(x') \mid \mb D_n)$ \\[0.5 ex]
& $\wht{f}_{n, \lambda}$ & KRR solution, $\argmin_{f \in \m H} \big[ n^{-1} \sum_{i=1}^n (Y_i - f(X_i))^2 + \lambda \norm{f}_{\m H}^2 \big]$ \\[0.5 ex]
& $\wt{K}$ & equivalent kernel, $\wt{K}(s,t) = \sum_{j=1}^{\infty} \nu_j \psi_j(s) \psi_j(t)$\\[0.5 ex]
& $\{\nu_j\}_{j=1}^{\infty}$ & eigenvalues of equivalent kernel, $\nu_j = \mu_j/(\mu_j + \lambda)$\\[0.5 ex]
& $F_{\lambda}$ & convolution with equivalent kernel, $F_{\lambda} g(t) = \int g(s) \wt{K}(s, t) ds$ \\[0.5 ex]
& $P_{\lambda}$ & $P_{\lambda} = \mbox{id} - F_{\lambda}$\\[0.5 ex]
& $S_{n, \lambda}$ & sample score function, $S_{n, \lambda} f = n^{-1} \sum_{i=1}^n (Y_i - f(X_i)) \wt{K}_{X_i} - P_{\lambda} f$ \\[0.5 ex]
& $S_{\lambda}$ & population score function, $S_{\lambda} f = F_{\lambda} f^\ast - f$ \\[0.5 ex]
& $U$ & second-order scaled bias approximation, $U(\cdot) = \sqrt{h/n} \sum_{i=1}^n w_i \wt{K}_{X_i}(\cdot)$\\[0.5 ex]
& $\wht{C}_n$ & covariance function of $U(\cdot)$, $\wht{C}_n(x, x') = \mb E U(x) U(x')$\\[0.5 ex]
& $\wht{W}$ & Gaussian approximation to $U$, $\wht{W} \sim \mbox{GP}(0, \wht{C}_n)$\\[0.5 ex]
& $\wt{Z}_n^M$ & rescaled sup-norm bias of posterior mean, $\wt{Z}_n^M = \sqrt{nh} \|\wht{f}_n - F_{\lambda} f^\ast\|_{\infty}$\\[0.5 ex]
& $Z_n^M$ & sup-norm of $U$, $Z_n^M = \|U\|_{\infty}$\\[0.5 ex]
& $\wht{Z}_n^M$ & sup-norm of $\wht{W}$, $\wht{Z}_n^M = \|\wht{W}\|_{\infty}$\\[0.5ex]
& $\wt{f}^B$ & $\wt{f}^B \sim \mbox{GP}(0, \wt{C}_n^B)$, law of centered GP posterior \\[0.5 ex]
& $\wht{C}_n^B$ & deterministic approximation to $nh \, \wt{C}_n^B$ \\[0.5 ex]
& $W^B$ & $W^B \sim \mbox{GP}(0, \wht{C}_n^B)$, approximation to law of $\sqrt{nh} \wt{f}^B$ \\[0.5 ex]
& $\wt{Z}_n^B$ & $\wt{Z}_n^B = \sqrt{nh} \, \|\wt{f}^B\|_{\infty}$, sup-norm of centered and scaled posterior GP \\[0.5 ex]
& $Z_n^B$ & $Z_n^B = \|W^B\|_{\infty}$, sup-norm of $W^B$ \\[0.5 ex]
\end{tabular}
\end{center}
\label{default}
\end{table}%


\subsection{Proof of Theorem~\ref{Thm:KRR_bound}}
\subsubsection{A First order bound for $\wht{f}_{n, \lambda}$}

Our analysis proceeds in two steps: we first prove a rough estimate for the error bound relative to the RKHS norm $\|\cdot\|$ using a sub-Gaussian type inequality, based on which we can control the covering entropy and further obtain a sharper sup-norm error bound using a Bernstein type inequality. We note that the sup-norm bound cannot be directly obtained by applying the Sobolev inequality that relates $\|\cdot\|$ with the sup norm.

\paragraph{A rough error bound relative to the $\|\cdot\|$ norm:} We make use of the following sub-Gaussian type concentration inequality under the $\|\cdot\|$ norm.
\begin{lemma}[Lemma 6.1 in \cite{yang2017non}]\label{Lemma:LD_bound}
Suppose $\alpha >1/2$. There exists some constant $C_K>0$ depending only on the kernel $K$, such that for any $t>0$,
\begin{align*}
\mb P \bigg( \sup_{g\in\m H}\bigg\{\Big\|\frac{1}{n}\,\sum_{i=1}^ng(X_i)\,\wt{K}_{X_i} - \mb E\big[g(X)\,\wt{K}_{X}\big]\Big\|\,\big/\,\|g\|\bigg\} \geq C_K\, n^{-\beta} \,t\bigg) \leq e^{-t^2},
\end{align*}
where the constant $\displaystyle \beta = \frac{(2\alpha-1)^2}{4\alpha\,(2\alpha +1)}$ is strictly positive.
\end{lemma}

Let $\m A_n$ denote an event such that under $\m A_n$, it holds for all $g\in\m H$,
\begin{align}\label{Eqn:EPbound}
\Big\|\frac{1}{n}\,\sum_{i=1}^ng(X_i)\,\wt{K}_{X_i} - \mb E\big[g(X)\,\wt{K}_{X}\big]\Big\| \leq C_K\, n^{-\beta/2}\,\|g\|.
\end{align}
According to Lemma~\ref{Lemma:LD_bound}, we may choose $\m A_n$ such that $\mb P(A_n)\geq 1- n^{-10}$ for any sufficiently large $n$.

Notice that we have the following two identities regarding the finite-sample level and the population level score functions,
\begin{align}\label{Eqn:FS_vs_P}
S_{n,\lambda}(\wht{f}_{n, \lambda}) = 0\quad\mbox{and} \quad S_{\lambda}(F_{\lambda} f^\ast) = 0.
\end{align}
Let $\Delta f=\wht{f}_{n, \lambda} - F_{\lambda} f^\ast$ denote the difference between the finite-sample minimizer and the population level minimizer of the KRR. From the identity~\eqref{eq:slam}, note also that $\Delta f = - S_{n, \lambda}(F_{\lambda} f^\ast)$. 
Now we will obtain a bound on $\Delta f$ by repeatedly applying inequality~\eqref{Eqn:EPbound}. 

By the definition of operators $S_{n,\lambda}$ and $S_\lambda$, we have
\begin{align}\label{Eqn:expre_1}
\big[S_{n,\lambda}(\wht{f}_{n, \lambda}) - S_\lambda(\wht{f}_{n, \lambda})\big] -\big[S_{n,\lambda}(F_{\lambda} f^\ast) - S_\lambda(F_{\lambda} f^\ast)\big] = \frac{1}{n}\,\sum_{i=1}^n\Delta f(X_i)\,\wt{K}_{X_i} - \mb E\big[\Delta f(X)\,\wt{K}_{X}\big].
\end{align}
Therefore, we obtain by applying inequality~\eqref{Eqn:EPbound} that under $\m A_n$, 
\begin{align*}
\Big|\big[S_{n,\lambda}(\wht{f}_{n, \lambda}) - S_\lambda(\wht{f}_{n, \lambda})\big] -\big[S_{n,\lambda}(F_{\lambda} f^\ast) - S_\lambda(F_{\lambda} f^\ast)\big] \Big| \leq C_K\, n^{-\beta/2}\,\|\Delta f\|.
\end{align*}
On the other hand, by using the identities in~\eqref{Eqn:FS_vs_P} along with the identity $\Delta f = - S_{n, \lambda}(F_{\lambda} f^\ast)$, we get
\begin{align}\label{Eqn:expre_2}
&\big[S_{n,\lambda}(\wht{f}_{n, \lambda}) - S_\lambda(\wht{f}_{n, \lambda})\big] -\big[S_{n,\lambda}(F_{\lambda} f^\ast) - S_\lambda(F_{\lambda} f^\ast)\big]  = \Delta f - S_{n,\lambda}(F_{\lambda} f^\ast). 
\end{align}
Let us know attempt to bound $\|S_{n,\lambda}(F_{\lambda} f^\ast)\|$. By applying the triangle inequality, we get
\begin{equation}\label{Eqn:S_n_bound}
\begin{aligned}
\norm{S_{n, \lambda}(F_{\lambda} f^\ast)} &= \big\|S_{n, \lambda}(F_{\lambda} f^\ast) - S_{\lambda}(F_{\lambda} f^\ast)  \big\| \\
 &\leq  \Big\|\frac{1}{n}\sum_{i=1}^n \{f^\ast - F_{\lambda} f^\ast\}(X_i)\, \wt{K}_{X_i} - \mb E\big[\{f^\ast - F_{\lambda} f^\ast\}(X)\,\wt{K}_{X}\big]  \Big\|    \\
 &\quad+
 \Big\|\frac{1}{n}\sum_{i=1}^n \{Y_i - f^\ast(X_i)\} \wt{K}_{X_i} \Big\|\\
 &\overset{(i)}{=} \Big\|\frac{1}{n}\sum_{i=1}^n P_\lambda f^\ast(X_i) \, \wt{K}_{X_i} - \mb E\big[P_\lambda f^\ast(X)\,\wt{K}_{X}\big]  \Big\|  +  
 \Big\|\frac{1}{n}\sum_{i=1}^n w_i\, \wt{K}_{X_i} \Big\|\\
 &:=  T_1 + T_2,
\end{aligned}
\end{equation}
where in step (i) we used the decomposition $\mbox{id}=P_\lambda + F_\lambda $, and recall that $w_i=Y_i-f^\ast(X_i)$ is the i.i.d.~$\m N(0,\,\sigma^2)$ noise.
Applying inequality~\eqref{Eqn:EPbound}, the first term $T_1$ can be bounded as
\begin{align*}
T_1 \leq C_K\, n^{-\beta/2}\,\|P_\lambda f^\ast\|.
\end{align*}
To bound the second term $T_2$, let $\Sigma=[\wt{K}(X_i,X_j)]_{1\le i,j\le n}$ and $w=(w_1,\ldots,w_n)^{\T}$. Then we have $T_2^2=n^2\, w^{\T}\Sigma w$.
By the Hanson-Wright inequality \citep{rudelson2013hanson}, we have
\[
P\Big[w^{\T}\Sigma w\geq\sigma^2\,\big( \mbox{tr}(\Sigma)+2\sqrt{\mbox{tr}(\Sigma^2)
\,t^2}+2\|\Sigma\|_F\, t^2\big)\Big]\leq e^{-t^2}, \quad \forall t>0,
\]
where $\|\cdot\|_F$ denotes the matrix Frobenius norm.
For some constant $D_K$ only depending on $K$, 
\begin{align*}
\mbox{tr}(\Sigma)&=\sum_{i=1}^n \wt{K}(X_i,X_i)\le D_K\, n\,h^{-1},\\
\mbox{tr}(\Sigma^2)&=\sum_{i,j=1}^n\wt{K}(X_i,X_j)^2\le D_K^2\,n^2\,h^{-2},\\
\|\Sigma\|_F&=\sqrt{\mbox{tr}(\Sigma^2)}\le D_K\,n\,h^{-1},
\end{align*}
we obtain
\[
P\left(
\Big\|\frac{1}{n}\sum_{i=1}^n w_i\,\wt{K}_{X_i}\Big\|\geq D_K\,\sigma\,(nh)^{-1/2}\,(1+2t+2t^2)
\right)\le e^{-t^2}.
\]
Therefore, there exists some event $\m B_n$ with $\mb P(B_n) \geq 1- n^{-10}$ such that on $\m B_n$ we have
\begin{align*}
T_2 \leq D_k\,\sigma\,(nh)^{-1/2}\,\log n,
\end{align*}
where we slightly abuse the notation by using $D_K$ to mean a different constant depending only on $K$.

By putting pieces together, we obtain that on $\m A_n\cap \m B_n$,
\begin{align*}
\|\Delta f\| &\leq \|\Delta f -  S_{n,\lambda}(F_{\lambda} f^\ast)\| + \| S_{n,\lambda}(F_{\lambda} f^\ast)\|\\
&\leq  C_K\, n^{-\beta/2}\,\|\Delta f\| +  C_K\, n^{-\beta/2}\,\|P_\lambda f^\ast\| + D_k\,\sigma\,(nh)^{-1/2}\,\log n.
\end{align*}
For $n$ sufficiently large such that $C_K\, n^{-\beta/2}\leq 1/2$, the above implies
\begin{align}\label{Eqn:First_Order}
\|\Delta f\| &\leq 2\, C_K\, n^{-\beta/2}\,\|P_\lambda f^\ast\| + 2\,D_k\,\sigma\,(nh)^{-1/2}\,\log n.
\end{align}
This further implies
\begin{align}
\|\wht{f}_{n, \lambda} -  f^\ast\| \leq \|\wht{f}_{n, \lambda} - F_{\lambda} f^\ast\| + \|f^\ast - F_{\lambda} f^\ast\| \leq 2\,\|P_\lambda f^\ast\| + 2\,D_k\,\sigma\,(nh)^{-1/2}\,\log n.
\end{align}

\paragraph{A sharp error bound relative to the sup norm:}
From inequality~\eqref{Eqn:First_Order} and the definition of $\|\cdot\|$, we obtain the following control on the Sobolev norm $\|\cdot\|_{\m H}$
\begin{align*}
\|\Delta f\|_{\m H} = \|\wht{f}_{n, \lambda} -  F_\lambda f^\ast\|_{\m H} \leq 2\lambda^{-1/2}\,\Big(\|P_\lambda f^\ast\| + D_k\,\sigma\,(nh)^{-1/2}\,\log n\Big):\, = \widetilde{A}_n.
\end{align*}
Moreover, we also have $\|\Delta f\|_{\infty} =\|\wht{f}_{n, \lambda} -  F_\lambda f^\ast\|_{\infty} \leq C\, h^{-1/2}\, \|\wht{f}_{n, \lambda} -  F_\lambda f^\ast\| \leq C'\, n$ for some sufficiently large constant $C'>0$.

For any $A_n >0$ and $B_n>0$, let $\m G_n = \{ g \in \m H : \norm{g}_{\infty} \le B_n,\,  \norm{g}_{\m H} \le A_n \}$.
We will make use of the following Bernstein type concentration inequality for suprema of empirical processes under the sup norm. A proof of Lemma \ref{Lemma:sup_norm_con} can be found in \S \ref{subsec:sup_norm_con}. 

\begin{lemma}\label{Lemma:sup_norm_con}
For any $\alpha >1/2$, it holds for some constant $C$ independent of $(n,h,A_n)$ that
\begin{align*}
\mb P \bigg[ &\sup_{t \in T; g \in \m G_n} \Big|\frac{1}{n} \sum_{i=1}^n g(X_i) \wt{K}(X_i, t) - E [ g(X_1) \wt{K}(X_1, t)]\Big|\\
&\  >   (nh)^{-1/2}\,\Big [\sqrt{\log n} + \sqrt{x}  + \sqrt{\log \frac{B_n}{ \norm{g}_{\infty}} }+ A_n^{1/(2\alpha)}\Big]\,  \norm{g}_{\infty}\\
 &\quad \ + (nh)^{-1}\, \Big[\log n + x +  \log \frac{B_n}{ \norm{g}_{\infty}} +A_n^{1/\alpha} \max\big(1,\, (n/h)^{1/(2\alpha) - 1/2}\big) \Big] \, \norm{g}_{\infty}\bigg] \leq 2e^{-x}, \quad x > 0.
\end{align*}
\end{lemma}

Let $\m B_n$ denote the event in this lemma with $A_n=\widetilde{A}_n$, $B_n=n$ and $t=c\,\log n$. Therefore, for sufficiently large $c$ we have $\mb P[\m B_n] \geq 1- n^{-10}$, and under event $\m B_n$ we have that for all $g\in \m H$ satisfying $\|g\|_\infty \leq n$ and $\|g\|_{\m H}\leq \widetilde{A}_n$,
\begin{equation}\label{Eqn:EP_sup_bound}
\begin{aligned}
&\Big\|\frac{1}{n} \sum_{i=1}^n g(X_i) \wt{K}_{X_i} - E [ g(X) \wt{K}_{X}]\Big\|_\infty \\
&  \leq C \, \Big[ \sqrt{\log \frac{n}{ \min(1,\, \norm{g}_{\infty})} } + (nh)^{-1/2}\,\max\big(1,\, (n/h)^{1/(2\alpha) - 1/2}\big) \, \widetilde{A}_n^{1/\alpha} \, \log \frac{n}{ \min(1,\, \norm{g}_{\infty})}\Big]\,\frac{\|g\|_\infty}{\sqrt{nh}} \\
& \leq C' \, \max\bigg( 1, \,n^{-1+1/(2\alpha)}\,h^{-1/(2\alpha)}\, \sqrt{\log \frac{n}{ \min(1,\, \norm{g}_{\infty})} }\bigg) \, \sqrt{\log \frac{n}{ \min(1,\, \norm{g}_{\infty})} }\ \frac{\|g\|_\infty}{\sqrt{nh}},
\end{aligned}
\end{equation}
for any $\alpha >1/2$ and $\lambda = h^{2\alpha} \asymp n^{-2\alpha/(2\alpha +1)}$.

Let us divide the proof into two cases:
\paragraph{If $\|\Delta f\|_\infty \leq \|P_\lambda f^\ast\|_\infty$:} The claimed bound follows since $\|\wht{f}_{n, \lambda} -  f^\ast\|_\infty \leq \|\Delta f\|_\infty + \|P_\lambda f^\ast\|_\infty$.

\paragraph{If $\|\Delta f\|_\infty \geq  \|P_\lambda f^\ast\|_\infty$:} In this case, we have $\log \frac{n}{ \min(1,\, \norm{\Delta f}_{\infty})} \asymp \log n$. 
Similar to the previous analysis, by applying~\eqref{Eqn:EP_sup_bound} with $g= \Delta f$, and comparing identities~\eqref{Eqn:expre_1} and \eqref{Eqn:expre_2}, we can get
\begin{align}\label{Eqn:errorbound_A}
\big\|\Delta f - S_{n,\lambda}(F_{\lambda} f^\ast)\big\|_\infty \leq C' \, \max\big( 1, \,n^{-1+1/(2\alpha)}\,h^{-1/(2\alpha)}\, \sqrt{\log n }\big) \, \sqrt{\log n }\ \frac{\|\Delta f\|_\infty}{\sqrt{nh}}.
\end{align}
In order to bound $\big\|S_{n,\lambda}(F_{\lambda} f^\ast)\big\|_\infty$, we will make use of the following concentration inequality for bounding $T_2$ relative to the sup norm.

\begin{lemma}\label{Lemma:empirical_kernel_bound}
 There exists some constant $C$ independent of $(n, h)$ such that for any $x>0$,
\begin{align}
 \bbP \bigg(\Big\|\frac{1}{n}\sum_{i=1}^n w_i\, \wt{K}_{X_i} \Big\|_{\infty}  \leq C\,\sigma\,\sqrt{\frac{1+x}{nh}} \bigg) \leq e^{-x}.  
 \end{align}
 \end{lemma}

By conducting a similar analysis as the steps in~\eqref{Eqn:S_n_bound} with $\|\cdot\|$ being replaced by $\|\cdot\|_\infty$, and using Lemma~\ref{Lemma:empirical_kernel_bound} and inequality~\eqref{Eqn:EP_sup_bound}, we can obtain that under intersection of event $\m B_n$ and the event in Lemma~\ref{Lemma:empirical_kernel_bound} with $x=c\,\log n$,
\begin{align}\label{Eqn:errorbound_B}
\big\|S_{n,\lambda}(F_{\lambda} f^\ast) \big\|_\infty \leq C\,\Big\{\max\big( 1, \,n^{-1+1/(2\alpha)}\,h^{-1/(2\alpha)}\, \sqrt{\log n }\big) \,\|P_\lambda f^\ast\|_\infty + \sigma \Big\}\,\sqrt{\frac{\log n}{nh}}.
\end{align}
Now combining inequalities~\eqref{Eqn:errorbound_A} and \eqref{Eqn:errorbound_B}, we get
\begin{align*}
&\big\|\Delta f\big\|_\infty \leq C' \,\big\|S_{n,\lambda}(F_{\lambda} f^\ast) \big\|_\infty 
 \leq C\,\Big\{\max\big( 1, \,n^{-1+1/(2\alpha)}\,h^{-1/(2\alpha)}\, \sqrt{\log n }\big) \,\|P_\lambda f^\ast\|_\infty + \sigma \Big\}\,\sqrt{\frac{\log n}{nh}},
\end{align*}
which implies the claimed bound as
\begin{align}\label{Eqn:First_Order_sup}
 \|\wht{f}_{n, \lambda} -  f^\ast\|_\infty \leq \|\Delta f\|_\infty + \|P_\lambda f^\ast\|_\infty \leq 2\|P_\lambda f^\ast\|_\infty  + C \, \sigma\, \sqrt{\frac{\log n}{nh}}.
\end{align}

\subsubsection{A second order bound for $\wht{f}_{n, \lambda}$}

Define $g_n := \wht{f}_{n, \lambda} -  f^\ast$.  By the definition of operators $S_{n,\lambda}$ and $S_\lambda$, we have
\begin{align}\label{Eqn:score_ind}
\big[S_{n,\lambda}(\wht{f}_{n, \lambda}) - S_\lambda(\wht{f}_{n, \lambda})\big] -\big[S_{n,\lambda}(f^\ast) - S_\lambda(f^\ast)\big] = \frac{1}{n}\,\sum_{i=1}^n g_n(X_i)\,\wt{K}_{X_i} - \mb E\big[g_n(X)\,\wt{K}_{X}\big].
\end{align}  
On the other hand, equation~\eqref{Eqn:FS_vs_P} leads to the identity
\begin{align*}
\big[S_{n,\lambda}(\wht{f}_{n, \lambda}) - S_\lambda(\wht{f}_{n, \lambda})\big] -\big[S_{n,\lambda}(f^\ast) - S_\lambda(f^\ast)\big] &= -S_\lambda(\wht{f}_{n, \lambda})+S_\lambda(f^\ast) - S_{n,\lambda}(f^\ast)\\
&\overset{(i)}{=} g_n - S_{n,\lambda}(f^\ast),
\end{align*}
where in step (i) we used \eqref{eq:slam}. 

Let us divide the proof into two cases:
\paragraph{If $\|g_n\|_\infty \leq h^{1/2}\, \|P_\lambda f^\ast\|_\infty$:} The claimed bound follows by using the fact that  $\Big\|\frac{1}{n}\,\sum_{i=1}^n g_n(X_i)\,\wt{K}_{X_i} - \mb E\big[g_n(X)\,\wt{K}_{X}\big]\Big\|_\infty \leq C\,h^{-1/2}\, \|g_n\|_\infty$.

\paragraph{If $\|\Delta f\|_\infty \geq  h^{1/2}\, \|P_\lambda f^\ast\|_\infty$:} 
In this case, we have $\log \frac{n}{ \min(1,\, \norm{g_n}_{\infty})} \asymp \log n$. Therefore, inequality~\eqref{Eqn:EP_sup_bound} and identity~\eqref{Eqn:score_ind} together imply that on $\m B_n$,
\begin{align*}
&\Big\|\big[S_{n,\lambda}(\wht{f}_{n, \lambda}) - S_\lambda(\wht{f}_{n, \lambda})\big] -\big[S_{n,\lambda}(f^\ast) - S_\lambda(f^\ast)\big] \Big\|_\infty \\
& \leq C' \,  \max\big( 1, \,n^{-1+1/(2\alpha)}\,h^{-1/(2\alpha)}\, \sqrt{\log n }\big) \, \sqrt{\log n }\ \frac{\|g_n\|_\infty}{\sqrt{nh}}.
\end{align*}

Combining the last two displays and the first order bound~\eqref{Eqn:First_Order_sup} for $\|g_n\|_\infty$, we obtain that for sufficiently large $n$, on event $\m A_n\cap\m B_n$,
\begin{equation}\label{Eqn:Second_Order}
\begin{aligned}
&\bigg\|\wht{f}_{n, \lambda} -  f^\ast -\bigg(\frac{1}{n}\,\sum_{i=1}^n w_i\,\wt{K}_{X_i} - P_\lambda f^\ast\bigg)\bigg\|_\infty \\
& \leq C' \, \max\big( 1, \,n^{-1+1/(2\alpha)}\,h^{-1/(2\alpha)}\, \sqrt{\log n }\big) \, \sqrt{\frac{\log n}{nh}} \,\Big(2\|P_\lambda f^\ast\|_\infty  + C \, \sigma\, \sqrt{\frac{\log n}{nh}} \Big).
\end{aligned}
\end{equation}
Moreover, since the randomness in $\m A_n\cap \m B_n$ is with respect to the noise $\{w_i:\,i=1,\ldots,n\}$ and the random design $\{X_i:\, i=1,\ldots,n\}$, this inequality also holds uniformly for all $f^\ast$ as well, meaning that it holds for all $f^\ast$ where $\wht{f}_{n, \lambda}$ is defined as a random function depending on $\{X_i:\, i=1,\ldots,n\}$ and $\{Y_i=f^\ast(X_i)+w_i:\,i=1,\ldots,n\}$.


\subsection{Proof of Corollary~\ref{cor:krr_supnorm}}
We prove the bounds on $\|P_{\lambda} f^\ast\|_{\infty}$ stated before the Corollary; the proof then is immediate upon optimally choosing $h$. 

First suppose $f^\ast \in \Theta_{\So}^{\alpha}(B)$. 
Bound $\|P_{\lambda} f^\ast\|_{\infty} \le \|P_{\lambda} f^\ast\| \, \sup_{x \in \m X} \| \wt{K}_x \|$. First, $\| \wt{K}_x \|^2 = \wt{K}(x, x) \lesssim \sum_{j=1}^{\infty} \nu_j \lesssim h^{-1}$ by \eqref{eq:nusum}. Hence, 
$\sup_{x \in \m X} \| \wt{K}_x \| \lesssim \lambda^{-1/(4 \alpha)}$. Next, 
$$
\| P_{\lambda} f^\ast \|^2 = \lambda \dotp{P_{\lambda} f^\ast}{f^\ast}_{\m H} = \lambda \sum_{j=1}^{\infty} \frac{\lambda}{\mu_j + \lambda} \ \frac{(f_j^\ast)^2}{\mu_j} \le \lambda \sum_{j=1}^{\infty} j^{2 \alpha} (f_j^\ast)^2. 
$$
This leads to $\|P_{\lambda} f^\ast \|_{\infty} \le \lambda^{1/2 - 1/(4 \alpha)} = h^{(\alpha - 1/2)}$ by using the definition of $\Theta_{\So}^{\alpha}(B)$. 

\

Now suppose $f^\ast \in  \Theta_{\H}^{\alpha}(B)$. Recall 
$$
P_{\lambda} f^\ast = \sum_{j=1}^{\infty} (1 - \nu_j) f_j^\ast \psi_j = \sum_{j=1}^{\infty} \frac{\lambda}{\mu_j + \lambda} \, f_j^\ast \psi_j.
$$
Bound 
$$
\|P_{\lambda} f^\ast\|_{\infty} \le \sum_{j=1}^{\infty} \frac{\lambda}{\mu_j + \lambda} f_j^\ast = \sqrt{\lambda} \sum_{j=1}^{\infty} \frac{ \sqrt{\lambda \mu_j}}{\mu_j + \lambda} \, \frac{f_j^\ast}{\sqrt{\mu_j}} \lesssim \sqrt{\lambda} \sum_{j=1}^{\infty} j^{\alpha} f_j^\ast \lesssim \sqrt{\lambda} = h^{\alpha},
$$
where we used the definition of $\Theta_{\H}^{\alpha}(B)$ and the AM-GM inequality to bound $\sqrt{\lambda \mu_j}/(\mu_j + \lambda) \le 1/2$.


\subsection{Proof of Theorem~\ref{thm:ratepw&sn}}
We prove the the first part.   
From the discussion in  \S  \ref{sec:Bayes_theory},  we obtain a high probability bound for the posterior expected point-wise loss.  With large probability
 \begin{align*}
\bbE[ |f(x^\ast) - f^\ast(x^\ast)|^2 \mid \mathbb{D}_n]  &= |\hat{f}_{n,\lambda}(x^\ast) - f^\ast(x^\ast)|^2 
+ \mbox{Var}[f(x^\ast) \mid  \mathbb{D}_n]  \\
 & \precsim \sigma^2\frac{\log n}{n\,h} + \|P_\lambda\,f^\ast\|_\infty^2  +   \|P_\lambda \,K_{x^\ast}\|^2_\infty.
\end{align*}
Since $\|P_\lambda \,K_{x^\ast}\|^2_\infty \precsim \sigma^2/(nh)$,  choosing $h$ optimally, we have \eqref{Eqn:fob1}.

To  prove \eqref{Eqn:fob2},  observe that 
 \begin{align*}
\bbE[ \|f -  f^*\|_{\infty}^2 \mid \mathbb{D}_n]  &\leq  2\|\hat{f}_{n,\lambda} - f^\ast\|_{\infty}^2 
+ 2\bbE \| f -  \hat{f}_{n,\lambda}\|_{\infty}^2.
\end{align*}
Since $f -  \hat{f}_{n,\lambda} \mid \mb D_n := Z \sim \mbox{GP}(0,  \wt{C}^B_n)$,  we use a version of Dudley's entropy integral (Theorem 3.2 of  \cite{pollard1989asymptotics}) for bounding moments of supremum of Gaussian processes. Define $\rho(s,t)  :=  \sqrt{\mbox{Var}\{ Z(s) - Z(t)\}}$ as the intrinsic psuedo-metric.   From the proof of Theorem 4.4., we have $\rho(s, t)  \leq  h^{-1} \abs{s - t}^{1/2}$.  
Hence the entropy of $[0, 1]^d$ with respect to $\rho$ is of the order $\log n$.   The width of the process $\{Z(t):  t \in [0, 1]^d\}$ with respect to $t_0 \in [0,1]^d$ is given by $\sup_t \rho(t, t_0)  \precsim \sigma/\sqrt{nh}$.  Thus from  Theorem 3.2 of  \cite{pollard1989asymptotics}
\begin{align*}
\bbE \| f -  \hat{f}_{n,\lambda}\|_{\infty}^2 \leq \frac{\sigma^2}{nh}  +   \frac{\sigma^2 \log  n}{nh}.
\end{align*}
Again choosing $h$ optimally, we have \eqref{Eqn:fob2}.


\subsection{Proof of Theorem~\ref{thm:point-wise_cov}}

Let us first find a normal approximation to the distribution of $\sqrt{nh}\,(\wht{f}_{n}(x) - F_\lambda f^\ast(x))$ for any fixed $x\in\m X$.
Using the Berry Esseen theorem, we obtain the following Kolmogorov distance bound between the sampling distribution of $\big\{\widehat{C}_n(x,x)\big\}^{-1/2}\sqrt{h/n}\,\sum_{i=1}^n w_i\,\wt{K}_{X_i}(x)$ and a standard normal random variable
\begin{align*}
\sup_{u\in\mb R}\bigg| \mb P\Big[ \big\{\widehat{C}_n(x,x)\big\}^{-1/2}\,\sqrt{h/n}\,\sum_{i=1}^n w_i\,\wt{K}_{X_i}(x) \leq u\Big] - \Phi(u)\bigg| \leq \frac{C}{\sqrt{nh}}, \quad\forall x\in\m X,
\end{align*}
where $\Phi$ denotes the cdf of a standard normal random variable, and we have used the fact that the third moment of $h^{1/2}w_i\,\wt{K}_{X_i}(x)$ is bounded by $h^{-1/2}$ up to a constant.
Since function $\Phi$ has a bounded derivative, we can combine the above with the second display in~\eqref{Eqn:mean_function_expansion} to obtain
\begin{align*}
\sup_{u\in\mb R}\bigg| \mb P\Big[ \big\{\widehat{C}_n(x,x)\big\}^{-1/2}\,\sqrt{nh}\,(\wht{f}_{n}(x) - F_\lambda f^\ast(x)) \leq u\Big] - \Phi(u)\bigg| \leq C\,\Big(\frac{1}{\sqrt{nh}} + \delta_n\Big), \quad\forall x\in\m X.
\end{align*}

Now let us turn to the posterior credible interval $\mbox{CI}_n(x;\,\beta)$. By combining the fact that $f(x)\sim \m N\big(\wht{f}_{n}(x), \,  \wt{C}^B_n(x,x)\big)$ and inequality~\eqref{Eqn:Cov_approx_bound}, we obtain
\begin{align*}
\sup_{u\in\mb R}\bigg| \mb P\Big[ \big\{\widehat{C}^B_n(x,\,x)\big\}^{-1/2} \sqrt{nh}\big(f(x) - \wht{f}_{n}(x)\big) \leq u\,\Big| \, \mb D_n\Big] - \Phi(u)\bigg| \leq C\,\gamma_n, \quad\forall x\in\m X.
\end{align*}
This implies that the half length of the credible interval $\mbox{CI}_n(x;\,\beta)$ satisfies
\begin{align*}
\bigg| l_n(x;\,\beta) - \sqrt{\frac{\widehat{C}^B_n(x,\,x)}{nh}} \, z_{(1+\beta)/2}\bigg| \leq C\,\sqrt{\frac{\widehat{C}^B_n(x,\,x)}{nh}}\, \gamma_n,
\end{align*}
where we used the fact the inverse function $\Phi^{-1}$ of $\Phi$ has a bounded derivative, and $z_u$ denotes the $u$th quantile of a standard normal distribution for $u\in(0,1)$.

We can express the frequentist coverage as 
\begin{align*}
&\mb P\big[f^\ast(x) \in  \mbox{CI}_n(x;\,\beta) \big] \\
= & \,\mb P\big[ - l_n(x;\,\beta)\leq  \wht{f}_{n}(x) - f^\ast(x)\leq  l_n(x;\,\beta)\big]\\
=& \,  \mb P\Big[   - \big\{\widehat{C}_n(x,x)\big\}^{-1/2}\,\sqrt{nh}\,l_n(x;\,\beta) + \big\{\widehat{C}_n(x,x)\big\}^{-1/2}\,\sqrt{nh}\, P_\lambda f^\ast(x) \\
&\qquad\qquad\qquad \leq  \big\{\widehat{C}_n(x,x)\big\}^{-1/2}\,\sqrt{nh}\,(\wht{f}_{n}(x) - F_\lambda f^\ast(x)) \\
&\qquad\qquad\qquad\qquad\qquad \leq  \big\{\widehat{C}_n(x,x)\big\}^{-1/2}\,\sqrt{nh}\,l_n(x;\,\beta) + \big\{\widehat{C}_n(x,x)\big\}^{-1/2}\,\sqrt{nh}\, P_\lambda f^\ast(x) \Big].
\end{align*}
Combining this with the previous approximation bounds, we can obtain
\begin{align*}
\bigg| &\mb P\big[f^\ast(x) \in  \mbox{CI}_n(x;\,\beta) \big] - \Phi\bigg(\sqrt{\frac{\widehat{C}^B_n(x,\,x)}{\widehat{C}_n(x,x)}}\, z_{(1+\beta)/2} + \big\{\widehat{C}_n(x,x)\big\}^{-1/2}\,\sqrt{nh}\, P_\lambda f^\ast(x)\bigg) \\
&\qquad + \Phi \bigg(-\sqrt{\frac{\widehat{C}^B_n(x,\,x)}{\widehat{C}_n(x,x)}}\, z_{(1+\beta)/2} + \big\{\widehat{C}_n(x,x)\big\}^{-1/2}\,\sqrt{nh}\, P_\lambda f^\ast(x)\bigg) \bigg| \leq C\, \Big( \frac{1}{\sqrt{nh}} + \gamma_n+\delta_n\Big),
\end{align*}
which implies the desired result.

\subsection{Proof of Corollary 3.1.} 
We define a sequence of integers $\{a_n\}$  to be of {\em constant separation}  if  $\{a_n\}$ is strictly increasing and $|a_n  - a_{n-1}| $ is constant for all $n\geq 2$. 
\\[1ex]
{\bf Under-smooth:}  
  For  any 
$x \in \mathcal{X}$,   there exists a constant $C_x$ depending on $x$, and a subsequence $k_j, j=1, \ldots, \infty$ with constant separation such that 
\begin{align*}
\widehat{C}_n(x,x)  \geq C_x h \sum_{j=1}^{\infty} \frac{1}{1+h^{2\alpha}/\mu_{k_j}}  = C_x.
\end{align*}
Also 
\begin{align*}
\sqrt{nh} \abs{P_\lambda f^*(x)}  \leq \lambda^{1/2} \sum_{j=1}^{\infty} \frac{j^{-\alpha_0}}{\lambda + j^{-2\alpha}}  \frac{|f_j^*|}{j^{-\alpha_0}}.
\end{align*}
For $2\alpha > \alpha_0$, the function $x \mapsto  x^{-\alpha_0}/(\lambda+ x^{-2\alpha})$ is maximized at $x = \{\alpha_0 \lambda / (2\alpha - \alpha_0) \}^{-1/(2\alpha)}$.   For $2\alpha \leq \alpha_0$, the function is monotonically decreasing.  In the first case $\sqrt{nh} \abs{P_\lambda f^*(x)} \leq C \lambda^{0.5(\alpha_0 / \alpha -1)}$ which goes to $0$.  In the later case,  
$\sqrt{nh} \abs{P_\lambda f^*(x)} \leq C \lambda^{1/2}$  which also goes to 0.  \\[1ex]
{\bf Smooth-match:} Let $(f^\ast_n)_{j}=C_\beta\,j^{-\alpha_0}$ for $j=j_n=\lfloor \lambda^{-1/(2\alpha)}\rfloor$, and $(f^\ast_n)_j=0$ for $j\neq j_n$, where $C_\beta$ is some tuning parameter to be determined later. It is obvious that $f^\ast\in \Theta_H^{\alpha_0}(B)$ for any $B\geq C_\beta$. Moreover, we have
\begin{align*}
\sqrt{nh} P_\lambda (f^\ast_n)(x)  = &\,\lambda^{1/2}\sum_{j=1}^{\infty} \frac{j^{-\alpha_0}}{\lambda + j^{-2\alpha}}  \frac{(f^\ast_n)}{j^{-\alpha_0}}\, \psi_j(x) = C_\beta \,\lambda^{1/2}\,\frac{\lfloor \lambda^{-1/(2\alpha)}\rfloor^{-\alpha}}{\lambda+\lfloor \lambda^{-1/(2\alpha)}\rfloor^{-2\alpha}}\,\psi_{j_n}(x) \\
:\,=&\, c_\lambda\,C_\beta\,\psi_{j_n}(x) \in [0.5\, C_\beta\,\psi_{j_n}(x),\, C_\beta\,\psi_{j_n}(x)].
\end{align*}
Let $\widetilde{b}_n$ be the solution of
\begin{align*}
\Phi(u_n(x;\,\beta) + \widetilde{b}_n) - \Phi((u_n(x;\,\beta) - \widetilde{b}_n) = \widetilde{\beta}.
\end{align*}
If we choose $C_\beta$ such that $\widetilde{b}_n=\widetilde{C}_n(x,x)^{-1/2}\,c_\lambda\,C_\beta\,\psi_{j_n}(x)$, then Theorem~\ref{thm:point-wise_cov} implies $\mb P\big[f_n^\ast(x) \in  \mbox{CI}_n(x;\,\beta) \big]\to \widetilde{\beta}$ as $n\to\infty$.
\\[1ex]
{\bf Over-smooth:}
Let $f^\ast_{bad,j}  =  (\mu_j^*)^{1/2} a_j$ where $ a_j =  j^{-(1+ \delta)}$ for   $\delta <  \alpha -  \alpha_0$.  Let 
$\beta:= 1+\delta+\alpha_0$.    There exists a subsequence $k_j$ with constant separation such that 
\small
\begin{align}\label{eq:os}
\sqrt{nh} P_\lambda f^\ast_{bad}(x) =  \lambda^{1/2} \sum_{j=1}^\infty \frac{(\mu_j^*)^{1/2}}{\lambda + \mu_j}  a_j  \psi_j(x) \gtrsim  \lambda^{1/2} \sum_{j=1}^\infty \frac{1}{\lambda k_j^{\beta} + k_j^{\beta - 2\alpha}}
&\asymp  \lambda^{\frac{1}{2} - \frac{1}{2\alpha} + \frac{\beta}{2\alpha} -1} \int \frac{dt}{t^{\beta} + t^{\beta - 2\alpha}} \nonumber \\ &\asymp  \lambda^{\frac{\delta + \alpha_0 - \alpha}{2\alpha}} \int \frac{dt}{t^{\beta} + t^{\beta - 2\alpha}}.  
\end{align}
\normalsize
Since  integral to the right of \eqref{eq:os} is finite and noting that $\delta + \alpha_0 - \alpha  < 0$, the conclusion follows immediately.

\subsection{Proof of Theorem~\ref{Thm:GP_mean_function}}
By applying Theorem 2.1 in \cite{chernozhukov2014gaussian} about a Gaussian approximation to the suprema of empirical processes (similar calculations as Corollary 2.2 therein), there exists a random variable $B_n$ that has the same distribution as $\widehat{Z}_n^M=\|\widehat{W}\|_\infty$, the suprema of the GP $\widehat{W}\sim \mbox{GP}(0,\,\widehat{C}_n)$, such that for any $\gamma\in(0,1)$,
 \begin{align*}
 \mb P\Big[\big| Z_n^M - B_n\big| \geq C_1\frac{b^{1/3}v^{2/3}\log^{2/3} n}{\gamma^{1/3}n^{1/6}}\Big] \leq C_2 \gamma,
 \end{align*}
where $b=h^{1/2} \sup_{x,x'}|\wt{K}(x,x')| \lesssim h^{-1/2}$, $v^2=\sup_x\mb E[h \,\wt{K}_X^2(x)]\lesssim 1$ and $C_1, C_2$ are constants independent of $n$. Combining this with Lemma 2.3 in \cite{chernozhukov2014gaussian} that converts this coupling bound to convergence in Kolmogorov distance, and the anti-concentration bound for GP in Corollary 2.1 in \citep{chernozhukov2014anti} that provides an upper bound to the pdf of $B_n$ or equivalently the pdf of $\widehat{Z}_n^M$, we can reach
\begin{align*}
\Big|\mb P\big[\widehat{Z}_n^M \leq t\big] - \mb P\big[Z_n^M \leq t\big] \Big| \leq C_3\,\bigg(\frac{\log^{2/3} n}{\gamma^{1/3}(nh)^{1/6}} + \gamma\bigg) \, A(|\widehat{W}|),
\end{align*}
where $A(|\widehat{W}|)=\mb E[Z_n^M] < \infty$ is a constant since the variance function $\widehat{C}_n(x,x)$ of the GP $\widehat{W}$ is uniformly bounded by some constant independent of $(n,h)$. Consequently, we obtain that for any $t\geq 0$,
\begin{align}
&\mb P\big[Z_n^M \leq t\big] - \mb P\big[\wt{Z}_n^M \leq t\big] \notag\\
& \overset{(i)}\leq \mb P\big[Z_n^M \leq t\big] - \mb P\big[\widehat{Z}_n^M \leq t +\sqrt{nh}\, \delta_n\big] + n^{-10} \notag\\
&\leq \Big|\mb P\big[Z_n^M \leq t\big] - \mb P\big[Z_n^M \leq t+\sqrt{nh}\, \delta_n\big]\Big| +  \Big|\mb P\big[Z_n^M \leq t+\sqrt{nh}\, \delta_n\big] - \mb P\big[\wt{Z}_n^M \leq t + \sqrt{nh}\, \delta_n\big]\Big| + n^{-10} \notag \\
&\overset{(ii)}{\leq} 4\,\sqrt{nh}\, \delta_n\, \big(A(|W|) + 1\big) + C_3\,\bigg(\frac{\log^{2/3} n}{\gamma^{1/3}(nh)^{1/6}} + \gamma\bigg) \, A(|W|) +  n^{-10}\notag\\
&\leq  C_4\, \bigg(\frac{\log^{2/3} n}{\gamma^{1/3}(nh)^{1/6}} + \gamma + \sqrt{nh}\, \gamma_n\, h^{2\alpha} +\gamma_n \sigma\,\sqrt{\log n}\,\bigg) + n^{-10}, \notag
\end{align}
where step (i) follows since the second display in~\eqref{Eqn:mean_function_expansion} implies $\big|\widehat{Z}_n^M - \wt{Z}_n^M\big| \leq \sqrt{nh}\,\delta_n$ with probability at least $1-n^{-10}$, and step (ii) follows by applying the anti-concentration bound in Corollary 2.1 in \citep{chernozhukov2014anti} for the suprema $Z_n^M$ of the GP $\widehat{W}$.
Similarly, we can obtain the same bound for $\mb P\big[\wt{Z}_n^M \leq t\big] -\mb P\big[Z_n^M \leq t\big]$, which combined with above implies that for any $t\geq 0$,
\begin{align}\label{Eqn:PosteriorModeLimit}
\Big|\mb P\big[Z_n^M \leq t\big] - \mb P\big[\wt{Z}_n^M\leq t\big] \Big| \lesssim \frac{\sqrt{\log n}}{(nh)^{1/8}},
\end{align}
by choosing $\gamma = \sqrt{\log n} / (nh)^{1/8}$.


\subsection{Proof of Theorem~\ref{Thm:GP_post}}\label{Sec:Proof_Thm_GP_post}
We need to apply the following Gaussian process comparison inequality, which is of independent interest. A proof is provided in \S~\ref{Sec:Proof_GP_comparison}.

\begin{theorem}[Comparison inequality for Gaussian processes]\label{Thm:GP_comparison}
Let $\|X\|_{\m F} = \sup_{f\in \m F} X_f$. Consider two centered Gaussian processes \mbox{GP}$(0,\, C^X)$ and \mbox{GP}$(0,\, C^Y)$ over an index set $\mathcal{F}$. Suppose $|C^X(f,\,g) - C^Y(f,\,g)| \leq \gamma$ for any $(f,g)\in \m F^2$. Let $\rho^X$ and $\rho^Y$ denote their respective intrinsic pseudometrics. If we let $Z_X:\,=\|X\|_{\m F}$ and $Z_X:\,=\|Y\|_{\m F}$, then for any $t\in\mb R$ and $\kappa\in(0,1)$,
\begin{align*}
\Big| \mb P\big[Z_X \leq t\big] - \mb P\big[Z_Y \leq t\big] \Big|& \leq 
\Big\{\nu_X+\nu_Y+\log N(\ve,\, \m F, \, \rho^X) + \log N(\ve,\, \m F, \, \rho^Y)\Big\}\,\gamma^{1/3}\\
&\qquad+\nu_X\,\delta_X+\nu_Y\,\delta_Y +2\ve\sqrt{2\log(4/\kappa)}+ \kappa,
\end{align*}
where
\begin{align*}
\nu_X&=K\, \sqrt{C^X(f,f)} + K \int_{0}^{\rho^X(\m F)}\sqrt{\log N(u,\, \m F,\, \rho^X)}\,du + 1\\
\nu_Y&=K\, \sqrt{C^Y(f,f)} + K \int_{0}^{\rho^Y(\m F)}\sqrt{\log N(u,\, \m F,\, \rho^Y)}\,du + 1\\
\delta_X&=K \int_0^\ve\sqrt{\log N(u,\, \m F,\, \rho^X)}\,du + \ve\,\sqrt{2\log(4/\kappa)}\\
\delta_Y&=K \int_0^\ve\sqrt{\log N(u,\, \m F,\, \rho^Y)}\,du + \ve\,\sqrt{2\log(4/\kappa)}.
\end{align*}
\end{theorem}

\noindent As a remark, for a GP over $T=[0,1]$, we can expand the index set to $T_E=T\cup [-2,-1]$ and define $X_t =-X_{t+2}$ for any $t\in[-2,-1]$. For this expanded index set, the conditions of Theorem~\ref{Thm:GP_comparison} still hold. In addition, this simple modification to the original GPs implies a Kolmogorov distance bound between their sup norms, that is, $\sup_{t\in T}|X_t|=\sup_{t\in T_E}X_t$ and $\sup_{t\in T}|Y_t|=\sup_{t\in T_E}Y_t$.

We will use this Gaussian process comparison inequality with $C^X=\wt{C}_n^B$ and $C^Y=\widehat{C}_n^B$ for proving the desired approximation error bound. Due to the sup norm error bound~\eqref{Eqn:Cov_approx_bound}, we can choose $\gamma \asymp \gamma_n$ in Theorem~\ref{Thm:GP_comparison}.
Before applying the Gaussian process comparison inequality, we also need study properties of intrinsic pseudometric $\widehat{\rho}_n(x, x')=\sqrt{nh}\,\sqrt{\mb E\big[(\wt{f}^B(x)-\wt{f}^B(x'))^2\big]}$ of the GP $\sqrt{nh}\,\wt{f}^B$ by comparing it with the Euclidean metric over $\m X=[0,1]$ (so that we can control the covering entropy of $\m X$ relative to the intrinsic pseudometric).
It is easy to see that we can express $\widehat{\rho}_n$ as
\begin{align}\label{Eqn:rho_n}
\sigma^{-2}\, \lambda\,h^{-1}\, \widehat{\rho}_n(x, x')^2=\big[g_{x,x'}(x)-g_{x,x'}(x')\big] - \big[\widehat{g}_{x,x'}(x)-\widehat{g}_{x,x'}(x')\big],
\end{align}
where $g_{x,x'}(\cdot) =  K(x, \cdot)-K(x',\cdot)$ and $\widehat{g}_{x,x'}$ is the solution of the following KRR with noiseless observations
\begin{align*}
\widehat{g}_{x,x'} = \argmin_{g \in \m H} \bigg\{ \frac{1}{n} \sum_{i=1}^n (Z^{x,x'}_i - g(X_i))^2+ \lambda \norm{g}_{\m H}^2 \bigg\}. 
\end{align*}
where $Z^{x,x'}_i = K(x,X_i)-K(x',X_i)  = g_{x,x'}(X_i)$. 
We can similarly apply the expansion~\eqref{Eqn:Second_Order} to obtain that with probability at least $1-n^{-10}$, 
\begin{align*}
\big\|\widehat{g}_{x, x'} - g_{x,x'} + P_\lambda \,g_{x,x'} \big\|_\infty \leq C\,\gamma_n\,\,\|P_\lambda \,g_{x,x'} \|_\infty.
\end{align*}
This implies that for all sufficiently large $n$, 
\begin{align*}
\big\|\widehat{g}_{x, x'} - g_{x,x'}\big\|_\infty &\leq 2 \,\|P_\lambda \,g_{x,x'}\|_\infty\\
&\leq 2\lambda \, \sup_{\tilde{x}} \Big|\sum_{j=1}^\infty \frac{\big\{\psi_j(x)-\psi_j(x')\big\}\, \psi_j(\tilde{x})}{1+h^{2\alpha}/\mu_j}\Big|\\
& \leq C\,\lambda\, L_\psi \,C_\psi\, |x-x'|\, \sum_{j=1}^\infty \frac{j}{1+(jh)^{2\alpha}}\\
&\lesssim \lambda\, h^{-2}\, |x-x'|.
\end{align*}
Combining this with identity~\eqref{Eqn:rho_n}, we obtain
\begin{align*}
&\sigma^{-2}\, \lambda\, h^{-1}\,\widehat{\rho}_n(x, x')^2 \lesssim 2\lambda \, h^{-2}\, |x-x'|, \\
\mbox{or}\qquad &\widehat{\rho}_n(x, x') \lesssim  h^{-1}\, |x-x'|^{1/2},\quad\forall (x,x')\in\m X^2.
\end{align*}
This implies that the covering entropy of $\m X=[0,1]$ relative to metric $\rho^X=\widehat{\rho}_n$ in Theorem~\ref{Thm:GP_comparison} satisfies $\log N(u,\m F,\rho^X) \lesssim \log(n/u)$ for any $u>0$. Similarly, it is straightforward to verify that the same bound holds for $\rho^Y$ as the intrinsic pseudometric associated with the population level GP $W^B$.  Moreover, it is easy to see that we can bound $\rho^X(\m F) = \sup_{x,x'} \widehat{\rho}_n(x, x')$ as
\begin{align*}
\sigma^{-2}\, \lambda\, h^{-1}\, \rho^X(\m F) &\leq 4 \sup_{x,x'}\,\big\|\widehat{g}_{x, x'} - g_{x,x'}\big\|_\infty\\
&\leq 4\lambda \, \sup_{\tilde{x}} \Big|\sum_{j=1}^\infty \frac{\big\{\psi_j(x)-\psi_j(x')\big\}\, \psi_j(\tilde{x})}{1+h^{2\alpha}/\mu_j}\Big|\\
&\leq C\,\lambda \, C_\psi^2 \,\sum_{j=1}^\infty \frac{1}{1+(jh)^{2\alpha}}\\
&\lesssim \lambda\, h^{-1},
\end{align*}
implying $\rho^X(\m F) \lesssim 1$. Similarly, it can be verified that $\rho^Y(\m F) \lesssim 1$.
Now we can apply the Gaussian comparison Theorem~\ref{Thm:GP_comparison} with $\gamma \asymp \gamma_n$, $\kappa\asymp n^{-1}$ and $\varepsilon\asymp \gamma_n^{1/3}$ to obtain that for any $t\geq 0$,
\begin{align}\label{Eqn:PosteriorLimit}
\Big|\mb P\big[\widetilde{Z}^B_n \leq t\, \big| \, \mb D_n\big] - \mb P\big[Z^B_n \leq t\big] \Big| \lesssim \gamma_n^{1/3} \,\log n.
\end{align}


\subsection{Proof of Theorem~\ref{thm:band_cov}}
By the definition of $r_n(\beta)$ as the $\beta$-th posterior quantile of $(nh)^{-1/2}\, \wt{Z}_n^B$, we obtain by choosing $t= \sqrt{nh}\,r_n(\beta)$ in inequality~\eqref{Eqn:PosteriorLimit_GP} that
\begin{align*}
\Big|\mb P\big[Z^B_n \leq \sqrt{nh}\,r_n(\beta)\big] - \mb P\big[Z^B_n \leq q_n^B(\beta)\big]\Big| \lesssim \gamma_n^{1/3}\log n,
\end{align*}
where we used the definition of $q_n^B(\beta)$ that implies $\mb P\big[Z^B_n \leq q_n^B(\beta)\big]=\beta$.
Now the anti-concentration inequality (Corollary 2.1 in \citep{chernozhukov2014anti}) for the supremum $Z_n^B$ of the GP $W^B$ implies that the pdf of $Z^B_n$ is bounded by some constant independent of $(n,h)$ because $A(|W^B|)=\mb E[Z_n^B] < \infty$ is a constant, we obtain that 
\begin{align*}
\Big| \sqrt{nh}\,r_n(\beta) -q_n^B(\beta)\Big| \lesssim \gamma_n^{1/3}\,\log n,
\end{align*}
implying the first desired bound comparing the quantiles.

Now we proceed to prove the second and third displayed inequalities in the theorem. By the definition, we have
\begin{align*}
\mb P\big[ f^\ast \in \mbox{CB}_n(\beta)\big] = \mb P\big[ \|\widehat{f}_n-f^\ast\|_\infty \leq r_n(\beta)\big] = \mb P\big[ \|\widehat{f}_n-F_\lambda f^\ast - P_\lambda f^\ast\|_\infty \leq r_n(\beta)\big]
\end{align*}
If the bias term $\|P_\lambda f^\ast\|_\infty$ satisfies $\sqrt{nh}\,\|P_\lambda f^\ast\|_\infty\to \infty$ as $n\to\infty$, then for arbitrarily large fixed constant $M>0$, as long as $n$ is sufficiently large, we always have the bound
\begin{align*}
\mb P\big[ f^\ast \in \mbox{CB}_n(\beta)\big]  \leq \mb P\big[ \|\widehat{f}_n-F_\lambda f^\ast\|_\infty \geq  \|P_\lambda f^\ast\|_\infty - r_n(\beta)\big] \leq \mb P\big[ \wt{Z}_n^M \geq  M\sqrt{nh}\, r_n(\beta)\big].
\end{align*}
Combining this with inequality~\eqref{Eqn:quantile_approx}, Theorem~\ref{Thm:GP_post}, and the Borell's inequality for GP, we obtain that for $n$ large enough,
\begin{align*}
\mb P\big[ f^\ast \in \mbox{CB}_n(\beta)\big]  \leq \mb P\big[ \wt{Z}_n^M \geq  M q_n^B(\beta)/2\big] \leq e^{-C\, M^2},
\end{align*}
for some constant $C$ independent of $(n,h)$. This implies the last desired bound~\eqref{Eqn:large_bias}.

On the other hand, if the bias term $\|P_\lambda f^\ast\|_\infty$ satisfies $\sqrt{nh}\,\|P_\lambda f^\ast\|_\infty\to 0$ as $n\to\infty$, then we have
\begin{align*}
\mb P\big[ f^\ast \in \mbox{CB}_n(\beta)\big]  \leq \mb P\big[  \|\widehat{f}_n-F_\lambda f^\ast\|_\infty \leq r_n(\beta) + \|P_\lambda f^\ast\|_\infty \big]=\mb P\big[ \wt{Z}_n^M \geq \sqrt{nh}\,(r_n(\beta) + \|P_\lambda f^\ast\|_\infty)\big].
\end{align*}
Now using Theorem~\ref{Thm:GP_post} and the anti-concentration inequality (Corollary 2.1 in \citep{chernozhukov2014anti}) for the supremum $Z_n^B$ of the GP $W^B$ lead to
\begin{align*}
\mb P\big[ f^\ast \in \mbox{CB}_n(\beta)\big]  - \mb P\big[ \widehat{Z}_n^M \geq \sqrt{nh}\, r_n(\beta) \big] \lesssim \sqrt{nh}\, \|P_\lambda f^\ast\|_\infty + \frac{\sqrt{\log n}}{(nh)^{1/8}}.
\end{align*}
A similar analysis leads to the same upper bound for $\mb P\big[ \widehat{Z}_n^M \geq \sqrt{nh}\, r_n(\beta) \big] - \mb P\big[ f^\ast \in \mbox{CB}_n(\beta)\big]$. These two together imply
\begin{align*}
\Big|\mb P\big[ f^\ast \in \mbox{CB}_n(\beta)\big]  - \mb P\big[ \widehat{Z}_n^M \geq \sqrt{nh}\, r_n(\beta) \big]\Big| \lesssim \sqrt{nh}\, \|P_\lambda f^\ast\|_\infty + \frac{\sqrt{\log n}}{(nh)^{1/8}}.
\end{align*}
Then the desired inequality~\eqref{Eqn:small_bias} is a direct consequence of the preceding display, inequality~\eqref{Eqn:quantile_approx} and the anti-concentration inequality (Corollary 2.1 in \citep{chernozhukov2014anti}) for the supremum $\widehat{Z}_n^M=\|\widehat{W}\|_\infty$ that implies a bounded pdf of $\widehat{Z}_n^M$ since $A(|\widehat{W}|)=\mb E[\widehat{Z}_n^M] < \infty$ is a constant.
Last but not least, since by construction the covariance function $\widehat{C}_n$ of the GP $\widehat{W}$ is uniformly strictly smaller than the covariance function $\widehat{C}_n^B$ of the GP $W^B$, we can conclude by a classical GP comparison inequality that $\|\widehat{W}\|_\infty$ is stochastically smaller than $\|\widehat{W}^B\|_\infty$. Since $q_n^B(\beta)$ is by definition the $\beta$-th quantile of $\|\widehat{W}^B\|_\infty$, we must have $\mb P\big[ \|\widehat{W}\|_\infty \leq q_n^B(\beta)\big] >\mb P\big[ \|\widehat{W}^B\|_\infty \leq q_n^B(\beta)\big] = \beta$.


\subsection{Proof of Theorem~\ref{Thm:GP_comparison}}\label{Sec:Proof_GP_comparison}

\
Let $\|X\|_{\m F_\ve}=\sup_{f,\, g\in\m F,\, \rho(f,\, g) \leq \ve } |X_f-X_g|$.

{\bf Step 1:} 
Fix $\ve>0$. Let $N = N(\ve,\, \m F, \, \rho^X) + N(\ve,\, \m F, \, \rho^Y)$, and $\{f_1,\ldots, f_N\}\subset \m F$ be the union of an $\ve$-net of $\m F$ relative to $\rho^X$ and an $\ve$-net of $\m F$ relative to $\rho^Y$. It is easy to see that $\{f_1,\ldots, f_N\}$ serves as an $\ve$-net of $\m F$ relative to both $\rho^X$ and $\rho^Y$.
To simplify the notation, we denote $X_{f_j}$ and $Y_{f_j}$ by $X_j$ and $Y_j$ respectively, for $j=1,2,\ldots,N$, and denote $C^X_{f_j,f_k}$ and $C^Y_{f_j,f_k}$ by $C^X_{jk}$ and $C^Y_{jk}$ for $j,k=1,\ldots,N$. 
Define
\begin{align*}
Z_X^\ve = \max_{1\leq j \leq N} X_j \quad\mbox{and}\quad Z_Y^\ve = \max_{1\leq j \leq N} Y_j,
\end{align*}
then $|Z_X^\ve - Z_X| \leq \|X\|_{\m F_\ve}$ and $|Z_Y^\ve - Z_Y| \leq \|Y\|_{\m F_\ve}$.
By applying Borell-Tsirelson-Ibragimov-Sudakov inequality~\cite[Proposition A.2.1]{van1996weak}, we obtain that with probability at least $1-\kappa/2$,
\begin{align*}
&\|X\|_{\m F_\ve} \leq \mb E[\|X\|_{\m F_\ve}] + \ve\,\sqrt{2\log(4/\kappa)},\quad\mbox{and}\\
&\|Y\|_{\m F_\ve} \leq \mb E[\|Y\|_{\m F_\ve}] + \ve\,\sqrt{2\log(4/\kappa)}.
\end{align*}
By the maximal inequality for GP~\cite[Corollary 2.2.8]{van1996weak}, we obtain
\begin{align*}
&E[\|X\|_{\m F_\ve}]\leq K \int_0^\ve\sqrt{\log N(u,\, \m F,\, \rho^X)}\,du,\quad\mbox{and}\\
&E[\|Y\|_{\m F_\ve}]\leq K \int_0^\ve\sqrt{\log N(u,\, \m F,\, \rho^Y)}\,du.
\end{align*}
Putting pieces together, we obtain that with probability at least $1-\kappa/2$,
\begin{equation}\label{Eqn:Dis_error}
\begin{aligned}
&|Z_X^\ve - Z_X|  \leq K \int_0^\ve\sqrt{\log N(u,\, \m F,\, \rho^X)}\,du + \ve\,\sqrt{2\log(4/\kappa)}=\delta_X,\quad\mbox{and}\\
&|Z_Y^\ve - Z_Y|  \leq K \int_0^\ve\sqrt{\log N(u,\, \m F,\, \rho^Y)}\,du + \ve\,\sqrt{2\log(4/\kappa)}=\delta_Y.
\end{aligned}
\end{equation}

{\bf Step 2:}
For any $t>0$, we approximate the non-smooth map $x\mapsto 1(\max_{1\leq j \leq N} x_j \leq t)$ by a smooth function. Following \cite{chernozhukov2013gaussian}, we first approximate the map $x\mapsto \max_{1\leq j \leq N} x_j$ by $F_\beta:\,\mb R^N\to \mb R$, defined as $F_\beta(x)=\beta^{-1}\log\big(\sum_{j=1}^Ne^{\beta x_j}\big)$. A straightforward calculation gives
\begin{align}\label{Eqn:Max_Ap_error}
\max_{1\leq j \leq N} x_j \leq F_\beta(x) \leq \max_{1\leq j \leq N} x_j + \beta^{-1}\log N.
\end{align}
Then we approximate the step function $x\mapsto 1(x\leq t)$ by $g(x) = v\big(\psi(x-t)\big)$ for some (large) $\psi>0$ and any smooth non-increasing function $v\in C^3(\mb R)$ satisfying $v(t)=1$ for $t\leq 0$, $v(t)\in[0,1]$ for $t\in[0,1]$, $v(t)=0$ for $t \geq 1$ and $\max\{\|v'\|_\infty,\,\|v''\|_\infty\|\}\leq C$ for some universal constant $C>0$.
It is straightforward to verify that
\begin{align}\label{Eqn:thre_Ap_error}
1(x\leq t) \leq g(x) \leq 1(x\leq t+\psi^{-1}), \quad \forall x\in\mb R,
\end{align}
and $\|g'\|_\infty\leq C\,\psi$, $\|g''\|_\infty\leq C\,\psi^2$.
Combining these two, we may approximate $x\mapsto 1(\max_{1\leq j \leq N} x_j \leq t)$ by the smooth function $f = g\circ F_\beta$. Our construction is closely related to the construction in \cite{chernozhukov2013gaussian} (they consider approximating the indicator function  $1_A(x)$ for any measurable subset $A$, therefore their construction is more complicated), and using their Lemma 4.3 we obtain
\begin{align}\label{Eqn:Hessian_bound}
\sum_{j=1}^N\sum_{k=1}^N |\partial_j\partial_k(g\circ F_\beta)(x)|\leq \|g''\|_\infty + 2\|g'\|_\infty \beta \leq C\,\big(\psi^2 + 2\psi \beta\big).
\end{align}

{\bf Step 3:} 
We need the following multivariate version of Stein's lemma, which can be proved by applying integration by parts.
\begin{lemma}\label{Lemma:Stein}
If $F:\,\mb R^N\to\mb R$ is a $C^1$ function with at most polynomial growth at infinity, and $W=(W_1,\ldots,W_n)$ is a centered Gaussian random vector, then for any $1\leq i\leq N$, 
\begin{align*}
\mb E\big[W_iF(W)\big]=\sum_{i=1}^N\mb E[W_iW_j]\,\mb E\big[\partial_i F(W)\big].
\end{align*}
\end{lemma}
Let $X=(X_1,\ldots,X_N)$ and $Y=(Y_1,\ldots,Y_N)$.
We consider the Slepian smart path interpolation $Z(t)=\sqrt{t}\, X+\sqrt{1-t}\, Y$ for $t\in[0,1]$, and let $\phi(t) = \mb E\big[g\circ F_\beta(Z_t)\big]$.
Then $\phi$ is differentiable, and
\begin{align*}
\phi'(t) = \mb E\bigg[\sum_{j=1}^N \partial_j(g\circ F_\beta)(Z_t) \,\Big(\frac{Y_j}{2\sqrt{t}} -\frac{X_i}{2\sqrt{1-t}}\Big)\bigg].
\end{align*}
For each $j$, Lemma~\ref{Lemma:Stein} implies
\begin{align*}
&\mb E\big[\partial_j(g\circ F_\beta)(Z_t) Y_j\big] =\sqrt{t}\, \sum_{k=1}^N C^Y_{jk}\,\mb E\big[\partial_{jk} (g\circ F_\beta)(Z_t) \big],\quad\mbox{and}\\
&\mb E\big[\partial_j(g\circ F_\beta)(Z_t) X_j\big] =\sqrt{1-t}\, \sum_{k=1}^N C^X_{jk}\,\mb E\big[\partial_{jk} (g\circ F_\beta)(Z_t) \big].
\end{align*}
Combining the last three displays, we obtain that for any $t\in[0,1]$,
\begin{align*}
\phi'(t) =\frac{1}{2}\,\sum_{j=1}^N\sum_{k=1}^N E\big[\partial_{jk} (g\circ F_\beta)(Z_t) \big]\,(C^X_{jk} - C^Y_{jk}).
\end{align*}
Combining this and inequality~\eqref{Eqn:Hessian_bound}, we finally reach
\begin{equation}\label{Eqn:App_boundA}
\begin{aligned}
&\big|\mb E\big[g\circ F_\beta(Y)\big] - \mb E\big[g\circ F_\beta(X)\big]\big| \\
\leq&\, \int_0^1|\phi'(t)|\,dt
\leq \frac{C}{2}\, (\psi^2+2\psi\beta)\, \big|C^X_{jk} - C^Y_{jk}\big| \leq \frac{C}{2}\, (\psi^2+2\psi\beta)\,\gamma,
\end{aligned}
\end{equation}
where in the last step we used the condition $\big|C^X_{jk} - C^Y_{jk}\big| \leq \gamma$.
\1

{\bf Step 4:}
Combining pieces together, we obtain that for any $t\in\mb R$ (notice that $g$ is non-increasing)
\begin{equation}\label{Eqn:key_diff}
\begin{aligned}
&\mb P[Z_X\leq t-\delta_X-\beta^{-1}\log N] -\mb P[Z_Y\leq t +\psi^{-1}+\delta_Y] &\\
\leq&\, \mb P[Z_X^\ve +\beta^{-1}\log N\leq t] - 
\mb P[Z_Y^\ve \leq t +\psi^{-1} ] + 2\ve\sqrt{2\log(4/\kappa)} + \kappa &\mbox{(by \eqref{Eqn:Dis_error})}\\
\leq&\,  \mb E\big[g(X+\beta^{-1}\log N)\big]-\mb E\big[g(Y)\big] + 2\ve\sqrt{2\log(4/\kappa)}+ \kappa   &\mbox{(by \eqref{Eqn:thre_Ap_error})}\\
\leq &\,  \mb E\big[g\circ F_\beta (X)\big]-\mb E\big[g\circ F_\beta(Y)\big] + 2\ve\sqrt{2\log(4/\kappa)}+ \kappa  &\mbox{(by \eqref{Eqn:Max_Ap_error})}\\
\leq &\, \frac{C}{2}\, (\psi^2+2\psi\beta)\,\gamma+ 2\ve\sqrt{2\log(4/\kappa)}+ \kappa.  &\mbox{(by \eqref{Eqn:App_boundA})}
\end{aligned}
\end{equation}
We have the anti-concentration bound for suprema of separable GP \cite[Theorem 2.1]{chernozhukov2014anti},
\begin{align*}
\mb P\big[|Z_X -t|\leq \epsilon\big] \leq 4\epsilon \,\mb E[Z_X] +4\epsilon,\quad\forall\epsilon>0,
\end{align*}
and the maximal inequality for GP \cite[Corollary 2.2.8]{van1996weak},
\begin{align*}
\mb E[Z_X] \leq K\, \sqrt{C^X(f,f)} + K \int_{0}^{\rho^X(\m F)}\sqrt{\log N(u,\, \m F,\, \rho^X)}\,du,
\end{align*}
where $\rho^X(\m F)=\sup_{f,g\in\m F}\rho^X(f,g)$ denotes the diameter of $\m F$ relative to $\rho^X$. Combining the two preceding displays, we obtain
\begin{align*}
\mb P\big[|Z_X -t|\leq \epsilon\big]  \leq 4\epsilon \,\Big\{K\, \sqrt{C^X(f,f)} + K \int_{0}^{\rho^X(\m F)}\sqrt{\log N(u,\, \m F,\, \rho^X)}\,du + 1\Big\}=\nu_X\,\epsilon.
\end{align*}
A similar anti-concentration bound holds for $Z_Y$ as
\begin{align*}
\mb P\big[|Z_Y -t|\leq \epsilon\big]  \leq 4\epsilon \,\Big\{K\, \sqrt{C^Y(f,f)} + K \int_{0}^{\rho^Y(\m F)}\sqrt{\log N(u,\, \m F,\, \rho^Y)}\,du + 1\Big\}=\nu_Y\,\epsilon.
\end{align*}
Putting pieces together, we finally reach that for any $t\in \mb R$,
\begin{align*}
&\mb P[Z_X\leq t] -\mb P[Z_Y\leq t]\\
\leq &\, \mb P[Z_X\leq t-\delta_X-(2\beta)^{-1}\log N-(2\psi)^{-1}] -\mb P[Z_Y\leq t +\delta_Y+(2\beta)^{-1}\log N+(2\psi)^{-1}]\\
\overset{(i)}{\leq} &\, \nu_X\,\big\{\delta_X+(2\beta)^{-1}\log N+(2\psi)^{-1}\big\}+ \nu_Y\,\big\{\delta_Y+(2\beta)^{-1}\log N+(2\psi)^{-1}\big\}\\
&\qquad\qquad\qquad\qquad\qquad\qquad\qquad\qquad\qquad +\frac{C}{2}\, (\psi^2+2\psi\beta)\,\gamma+ 2\ve\sqrt{2\log(4/\kappa)}+ \kappa,
\end{align*}
where in step (i) we applied inequality~\eqref{Eqn:key_diff} with $t=t+(2\beta)^{-1}\log N-(2\psi)^{-1}$.
Now, we can take $\beta = \psi \log N$ and $\psi = \gamma^{-1/3}$ to obtain
\begin{align*}
\mb P[Z_X\leq t] -\mb P[Z_Y\leq t] \leq \nu_X\,\delta_X+\nu_Y\,\delta_Y + 
(\nu_X+\nu_Y+\log N)\,\gamma^{1/3}+2\ve\sqrt{2\log(4/\kappa)}+ \kappa.
\end{align*}

\subsection{Proof of Lemma~\ref{Lemma:empirical_kernel_bound}}

To begin with, we collect a version of the classical Bernstein inequality and a resulting expectation bound. \\
{\bf Bernstein's inequality:} Let $X_1, \ldots, X_n$ be independent random variables. Let $\nu$ and $c$ be positive numbers such that  $\sum_{i=1}^n \mb E X_i^2 \le \nu$, and $\sum_{i=1}^n \mb E |X_i|^q \le 2^{-1} q! ~ \nu c^{q-2}$ for $q \ge 3$. Let $S = \sum_{i=1}^n (X_i - \mb E X_i)$. 
Then, for any $\lambda \in (0, 1/c)$, 
$$
\mb E e^{\lambda S} \le \exp \bigg[\frac{\nu \lambda^2}{2(1 - c \lambda)}\bigg], 
$$
which in particular implies, for any $x > 0$, 
\begin{align}\label{eq:bern}
\mb P[ |S| \ge \sqrt{2 \nu x} + c x] \le 2 e^{-x}. 
\end{align}
Further, if $S$ a random variable satisfying \eqref{eq:bern}, then 
Then, for any $q \ge 1$, 
$$
\mb E[S^{2 q}] \le q! (8 \nu)^{q} + (2 q!) (4c)^{2 q}. 
$$
In particular, $\mb E (S^2) \lesssim C_1 \nu + C_2 c^2$, where $C_1, C_2$ are absolute constants. 

Our proof makes use of the following tail bound for supremum of empirical processes with sub-exponential increments from \cite{baraud2010bernstein}. \\
{\bf Bernstein-type inequality for suprema of random processes (Theorem 2.1, \cite{baraud2010bernstein})} Let $(X_t)_{t \in T}$ be a centered family with $T \subset \mb R^D$ for some finite $D$. Fix some $t_0 \in T$ and let $\bar{Z} = \sup_t |X_t - X_{t_0}|$. Consider norms $d(s, t) =d(s-t)$ and $\delta(s, t) = \delta(s-t)$ on $\mb R^D$, and assume there exist $v, b > 0, c \ge 0$ such that 
$$
T \subset \{ t \in \mb R^D : d(t, t_0) \le v, c \delta(t, t_0) \le b \}.
$$
Further, assume that for all $s \ne t \in T$, 
\begin{align}\label{eq:increm_bern}
\mb E \big[ e^{\lambda (X_t - X_s)} \big] \le \exp \bigg[ \frac{\lambda^2 d^2(t, s)}{2 (1 - \lambda c \delta(t, s))} \bigg],  \ \forall \ \lambda \in \bigg[0, \frac{1}{c \delta(t, s)}\bigg]. 
\end{align}
Then, with $C = 18$, 
\begin{align}\label{eq:EPbern}
\mb P\big[ \bar{Z} \ge C \big( \sqrt{ v^2 (1 + x) } + b (1 + x) \big) \big] \le 2 e^{-x}, \ \forall \ x > 0. 
\end{align}

Without loss of generality, we assume $\sigma=1$.
Let $U_t = \sum_{i=1}^n w_i \wt{K}(X_i, t)$, where $w_i \sim N(0, 1)$, $X_i \sim U(0, 1)$, and $w_i $'s are independent of $X_i$'s. It suffices to find a tail bound for $\norm{U}_{\infty} = \sup_{t \in [0, 1]} |U_t|$. For $t, s \in [0, 1]$, we write $U_t - U_s = \sum_{i=1}^n w_i \Delta_i^{t,s}$, where $\Delta_i^{t,s} = \wt{K}(X_i, t) - \wt{K}(X_i, s) = \sum_{j=1}^{\infty} \nu_j \{ \psi_j(t) - \psi_j(s) \} \psi_j(X_i)$. We suppress the dependence on $s$ and $t$, and write $\Delta_i$ subsequently. 

We now proceed to the proof. Recall $U_t - U_s = \sum_{i=1}^n w_i \Delta_i$. We first define the two norms $d$ and $\delta$. Let $d(t, s) = [ E(U_t - U_s)^2 ]^{1/2}$ be the intrinsic semi-metric. For any $\kappa \in (0, \alpha \wedge 1)$, let $\delta(t, s) = |t - s|^{\kappa}$. Now we estimate the quantities $v, b$ and $c$ appearing in Bernstein's inequality~\eqref{eq:EPbern}. First, we notice 
\begin{align*}
d^2(t, s) = E \sum_{i=1}^n (w_i \Delta_i)^2 = \sum_{i=1}^n E \Delta_i^2. 
\end{align*}
Under Assumption B, we can obtain by using orthonormality of the eigenfunctions that
\begin{align*}
\mb E \Delta_i^2 = \sum_{j=1}^{\infty} \nu_j^2 |\psi_j(t) - \psi_j(s)|^2 \le 2 C_{\psi}^2 \sum_{j=1}^{\infty} \nu_j^2 \le C h^{-1}. 
\end{align*}
Combining the two previous displays, we conclude 
\begin{align*}
\sup_{t, s} d(t, s) \le \sqrt{n h^{-1}}. 
\end{align*}
In order to verify condition~\eqref{eq:increm_bern} through characterizing the growth of the moments, we need to bound $E \sum_{i=1}^n |w_i \Delta_i|^q$ for $q \ge 3$. To that end, we first bound $|\Delta_i|$ as follows. 
For any $i = 1, \ldots, n$ and any fixed $\kappa\in(0,1)$, we use Assumption B to get
\begin{align*}
|\Delta_i| & \le \sum_{j=1}^{\infty} \nu_j |\psi_j(t) - \psi_j(s)| |\psi_j(X_i)|  \\
& \le C_{\psi} \sum_{j=1}^{\infty} \nu_j |\psi_j(t) - \psi_j(s)|^{\kappa} [2 C_{\psi}^{1 - \kappa}]  \\
& \le 2 C_{\psi}^{2 - \kappa} L_{\psi}^{\kappa} ~ |t - s|^{\kappa} \sum_{j=1}^{\infty} j^{\kappa} \nu_j \le C h^{-(1 + \kappa)} \delta(t, s).
\end{align*}
This implies that for any $q \ge 3$, 
\begin{align*}
\mb E \sum_{i=1}^n |w_i \Delta_i|^q &\le \sum_{i=1}^n q^{q/2} \mb E \big[ |\Delta_i|^{q-2} \Delta_i^2] \le 2^{-1} q! \ d^2(t, s) \ \big[ h^{-(1+\kappa)} \delta(t, s) \big]^{q-2},
\end{align*}
where we used the fact that $\mb E |w_i|^q \lesssim q^{q/2}$ and $q! \ge (q/e)^q$. For the $|\Delta_i|^{q-2}$ term, we used the global bound on $|\Delta_i|$, and finally used $\sum_{i=1}^n \mb E \Delta_i^2 = d^2(t, s)$. 

Let $c = h^{-(1 + \kappa)}$. Then using the classical Bernstein's inequality described in the beginning of this subsection, we obtain 
$$
E e^{\lambda (U_t - U_s)} \le \exp \bigg[ \frac{\lambda^2 d^2(t, s)}{2 (1 - \lambda c \delta(t, s))} \bigg], \ \forall \lambda \in \bigg[0, \frac{1}{c \delta(t, s)} \bigg]. 
$$
Thus the sub-exponential increment condition~\eqref{eq:increm_bern} is satisfied, and quantities $v^2$ and $b$ in Bernstein's inequality~\eqref{eq:EPbern} are $n h^{-1}$ and $h^{-(1+\kappa)}$, respectively. Thus, 
$$
\mb P\big[ \|U\|_\infty \ge 18 \big( \sqrt{ n h^{-1} (1 + x) } + h^{-(1+\kappa)} (1 + x) \big) \big] \le 2 e^{-x}, \ \forall \ x > 0,
$$
which implies the claimed concentration inequality by dividing both sides inside the probability by $n$.

\subsection{Proof of Lemma~\ref{Lemma:sup_norm_con}}\label{subsec:sup_norm_con}

Define 
\begin{align*}
U(t, g) = \frac{1}{n} \sum_{i=1}^n g(X_i) \wt{K}(X_i, t) - E [ g(X_1) \wt{K}(X_1, t)]. 
\end{align*}

To prove the result, we will apply the peeling technique and the following lemma. A proof of the following lemma is deferred to the next subsection. In this subsection, the meaning of constant $C$ can be change from line to line.

\begin{lemma}\label{Lemma:fixed_rho}
Let $\m G_n(\rho) = \{ g \in \m H : \norm{g}_{\infty} \le \rho,\,  \norm{g}_{\m H} \le A_n \}$. Then
\begin{align*}
P \bigg( \sup_{t \in T; g \in \m G_n} |U(t, g)| >  H_n+ (nh)^{-1/2}\rho \sqrt{x}  +  (nh)^{-1}\rho\, x \bigg) \leq e^{-x}, \quad x > 0, 
\end{align*}
where $H_n =  (nh)^{-1/2} \rho\, [1+ \sqrt{\log n} + A_n^{1/(2\alpha)}] 
+ (nh)^{-1} \rho\, [1+ \log n + A_n^{1/\alpha} \max\big(1,\, (n/h)^{1/(2\alpha) - 1/2}\big)] $ for any $\alpha > 1/2$. 
\end{lemma}

Now we proceed to the proof. To apply the peeling technique, we decompose the range $(0, B_n]$ into $\bigcup_{j=0}^\infty (2^{-j-1} B_n,\, 2^{-j} B_n]$. For any $g\in \m G_n$, there always exists some $j_0\geq 0$ such that
$\|g\|_\infty \in (2^{-j_0-1} B_n,\, 2^{-j_0} B_n]$, which implies $j_0 \leq c_1 \sqrt{\log \frac{B_n}{ \norm{g}_{\infty}} }$ for some constant $c_1 >0$. Using Lemma~\ref{Lemma:fixed_rho} with $\rho = 2^{-j_0} B_n$ and $x=x+j_0$, we obtain
\begin{align*}
P \bigg(& \sup_{t \in T; g \in \m G_n,\, \|g\|_\infty \in (2^{-j_0-1} B_n,\, 2^{-j_0} B_n]} |U(t, g)| > C\, (nh)^{-1/2}\,\Big [\sqrt{\log n} + \sqrt{x} + \sqrt{\log \frac{B_n}{ \norm{g}_{\infty}} } + A_n^{1/(2\alpha)}\Big]\,  \norm{g}_{\infty} \\
 &\qquad \qquad + C\, (nh)^{-1}\, \Big[\log n + x +\log \frac{B_n}{ \norm{g}_{\infty}} + A_n^{1/\alpha} \max\big(1,\, (n/h)^{1/(2\alpha) - 1/2}\big) \Big] \, \norm{g}_{\infty} \bigg) \leq \frac{1}{2^{j_0}}\, e^{-x},
\end{align*}
where we used the inequality that $\sqrt{a+b} \leq \sqrt{a} + \sqrt{b}$ and the fact that $\|g\|_\infty >2^{-j_0-1} B_n = \rho/2$.
Now by combining this with a union bound over $j_0$, we obtain
\begin{align*}
P \bigg(& \sup_{t \in T; g \in \m G_n} |U(t, g)| > C\, (nh)^{-1/2}\,\Big [\sqrt{\log n} + \sqrt{x} + \sqrt{\log \frac{B_n}{ \norm{g}_{\infty}} } + A_n^{1/(2\alpha)}\Big]\,  \norm{g}_{\infty} \\
 &\qquad \qquad + C\, (nh)^{-1}\, \Big[\log n + x +\log \frac{B_n}{ \norm{g}_{\infty}} + A_n^{1/\alpha} \max\big(1,\, (n/h)^{1/(2\alpha) - 1/2}\big) \Big] \, \norm{g}_{\infty} \bigg)\\
 \leq \sum_{j_0=0}^\infty P \bigg(& \sup_{t \in T; g \in \m G_n,\, \|g\|_\infty \in (2^{-j_0-1} B_n,\, 2^{-j_0} B_n]} |U(t, g)| > C\, (nh)^{-1/2}\,\Big [\sqrt{\log n} + \sqrt{x} + \sqrt{\log \frac{B_n}{ \norm{g}_{\infty}} } + A_n^{1/(2\alpha)}\Big]\,  \norm{g}_{\infty} \\
 &\qquad \qquad + C\, (nh)^{-1}\, \Big[\log n + x +\log \frac{B_n}{ \norm{g}_{\infty}} + A_n^{1/\alpha} \max\big(1,\, (n/h)^{1/(2\alpha) - 1/2}\big) \Big] \, \norm{g}_{\infty} \bigg) \\
 &  \leq \sum_{j_0=0}^\infty \frac{1}{2^{j_0}}\, e^{-x} =2 e^{-x},
\end{align*}
which is the desired result.

\subsection{Proof of Lemma~\ref{Lemma:fixed_rho}}
Let $\eta = (t, g) \in \m T \otimes \m H$. The proof applies an improved version of Bernstein's inequality~\eqref{eq:EPbern} over the product space $\m T\otimes \m H$ by truncating the chaining in the proof of Theorem 5.1 in \cite{baraud2010bernstein} at a finite level.
We will use the following simple inequality multiple times: 
\begin{align*}
\mb E | g(X_1) \wt{K}(X_1, t) | \le C \norm{g}_{\infty} h^{-1/2}, \quad\mbox{for any $g \in \m H$.}
\end{align*}
In fact, this follows since  $[ \mb E | g(X_1) \wt{K}(X_1, t) | ]^2 \le \mb E [ g^2(X_1) \wt{K}^2(X_1, t)]  \le \norm{g}_{\infty}^2 \mb E \wt{K}^2(X_1, t)  \le C \norm{g}_{\infty}^2 h^{-1}$.

First, set $d^2(\eta, \eta') =\mb E [ U(t, g) - U(t', g') ]^2$ to be the intrinsic semi-metric in Bernstein's inequality~\eqref{eq:EPbern} as in the proof of Lemma~\ref{Lemma:empirical_kernel_bound}. Let's try to get a global bound for $d$ first for finding $v$,
\begin{align*}
d^2(\eta, \eta') 
& = \mbox{Var} [ U(t, g) - U(t', g') ] \\
& \le 2 \big[ \mbox{Var} U(t, g) + \mbox{Var} U(t', g') \big] \\
& \le \frac{2}{n} E \big[ g^2(X_1) \wt{K}^2(X_1, t) + g'^2(X_1) \wt{K}^2(X_1, t') \big] \\
& \le \frac{2 h^{-1} ( \norm{g}_{\infty}^2 + \norm{g'}_{\infty}^2)}{n}. 
\end{align*}
As before, we need verify condition~\eqref{eq:increm_bern} by establishing a Bernstein type inequality for the difference $U(t, g) - U(t', g')$. Write $U(t, g) - U(t', g') = \sum_{i=1}^n W_i$, where 
\begin{align*}
W_i = n^{-1} \bigg[ \big\{ g(X_i) \wt{K}(X_i, t) - g'(X_i) \wt{K}(X_i, t') \big\} - \mb E \big\{ g(X_1) \wt{K}(X_1, t) - g'(X_1) \wt{K}(X_1, t') \big\} \bigg]. 
\end{align*} 
Clearly, $\sum_{i=1}^n \mb E W_i^2 = d^2(\eta, \eta')$. We first obtain a bound on $|W_i|$ in order to find quantity $c$ and metric $\delta$ in~\eqref{eq:increm_bern}. To simplify the notation, we write $V_{g, t}(X_i) = g(X_i) \wt{K}(X_i, t) - \mb E g(X_1) \wt{K}(X_1, t)$. 
\begin{align*}
n |W_i| = |V_{g,t}(X_i) - V_{g',t'}(X_i)| \le |V_{g,t}(X_i) - V_{g,t'}(X_i)| + |V_{g,t'}(X_i) - V_{g',t'}(X_i)|. 
\end{align*}
We bound the two terms on the right hand side of the above display separately. Under Assumption B,  we have for any fixed $\kappa\in(0,1)$,
\begin{align*}
|V_{g,t}(X_i) - V_{g,t'}(X_i)| 
&\le \abs{ g(X_i) \{ \wt{K}(X_i, t) - \wt{K}(X_i, t') \} } + \mb E \abs{ g(X_1) \{ \wt{K}(X_1, t) - \wt{K}(X_1, t') \} } \\
 & \le \norm{g}_{\infty} h^{-(1+\kappa)} |t - t'|^{\kappa} + \bigg[ \norm{g}_{\infty}^2 h^{-(1+2\kappa)} |t - t'|^{2 \kappa} \bigg]^{1/2} \\
 & \lesssim \norm{g}_{\infty} h^{-(1+\kappa)} |t -t'|^{\kappa}. 
\end{align*}
In addition, under the same assumption we have
\begin{align*}
|V_{g,t'}(X_i) - V_{g',t'}(X_i)| 
& \le \abs{ (g - g')(X_i) \wt{K}(X_i, t') } + \mb E \abs{ (g - g')(X_1) \wt{K}(X_1, t') \} }\\
& \le \norm{g - g'}_{\infty} h^{-1} + \bigg[ \norm{g - g'}_{\infty}^2 h^{-1} \bigg]^{1/2} \\
& \lesssim \norm{ g - g'}_{\infty} h^{-1}. 
\end{align*}
Putting these together, we obtain
\begin{align*}
|W_i| \lesssim (n h)^{-1} \big[ \norm{g - g'}_{\infty}  + h^{-\kappa} \norm{g}_{\infty} |t - s|^{\kappa}  \big].
\end{align*}
This means we can choose 
\begin{align}
\delta(\eta, \eta') = \big[ \norm{g - g'}_{\infty}  + h^{-\kappa} \norm{g}_{\infty} |t - s|^{\kappa}  \big],
\end{align}
and $c=(nh)^{-1}$.
Under this choice of $c$ and $\delta$, we can easily verify that for $q \ge 3$, 
\begin{align*}
\mb E \sum_{i=1}^n |W_i|^q \le [(nh)^{-1} \delta(\eta, \eta')]^{q-2} d^2(\eta, \eta'),
\end{align*}
which implies the sub-exponential increment condition~\eqref{eq:increm_bern}.
 To calculate covering numbers in Theorem 5.1 of \cite{baraud2010bernstein}, we need to bound the $d$ metric by the Euclidean metric. To that end, 
\begin{align*}
d^2(\eta,\eta') 
& \le 2 \big[ \mbox{Var} \{ U(t, g) - U(t',g) \} + \mbox{Var} \{(U(t',g) - U(t',g') \} \big] \\
&\le \frac{2}{n} \big[ h^{-(1+2\kappa)} \norm{g}_{\infty}^2 |t - s|^{2 \kappa} + h^{-1} \norm{g - g'}_{\infty}^2 \big]. 
\end{align*}
Therefore, on the product space $T \otimes \m G$, we can choose $v = (nh)^{-1/2}\rho, c = (nh)^{-1}$ and $b = (nh)^{-1}\rho$ in Theorem 5.1 of \cite{baraud2010bernstein}. Moreover, using their notations, we have the following bound on the level $k$ set $\m A_k$ for telescoping,
$$
|\m A_k| \le \max\big\{  N(\m T \otimes \m G, d, 2^{-k} v),  N(\m T \otimes \m G, c \delta, 2^{-k} b) \big\},
$$
where $N(\m F, d, u)$ denotes the $u$-covering number of space $\m F$ relative to metric $d$.
In our case, since by definition $\m G$ is the $\alpha$th order Sobolev space with radius $A_n$, we have
\begin{align*}
 N(\m T \otimes \m G, d, 2^{-k} v)  \leq \exp\{(2^k A_n)^{1/\alpha}\} (2^k / h^\kappa)^{1/\alpha}, \quad 
 N(\m T \otimes \m G, c\delta, 2^{-k} b)  \leq \exp\{(2^k A_n)^{1/\alpha}\} (2^k / h^\kappa)^{1/\alpha}.
\end{align*}

Now we describe an improved version of Theorem 5.1 in \cite{baraud2010bernstein} that we will use. We still apply the chaining technique, but applying an improved telescoping identity by truncating at some finite level $k_0$ (please refer to the first display in the  proof of Theorem 5.1 in \S 5.1 in \cite{baraud2010bernstein}, where we have adopted our notation by identifying $X_{\cdot}$ with $U(\cdot)$ and $t$ with $\eta=(t,g)$),
\begin{align*}
U(\eta) - U(\eta_0) = \sum_{k=0}^{k_0 - 1}\Big[U(\pi_{k+1}(\eta)) - U(\pi_{k}(\eta)) \Big]+ \big[ U(\eta) - U(\pi_{k_0}(\eta))\big],\quad \forall \eta\in \m T \otimes \m H,
\end{align*}
where $\pi_k$ maps any point $\eta\in\m T \otimes \m H$ to some point $\pi_k(\eta)$ in $\m A_k$ such that $d\big(\eta, \pi_k(\eta)\big) \leq 2^{-k}v$ and $c\delta\big(\eta, \pi_k(\eta)\big)\leq 2^{-k} b$, and satisfies $\pi_0(\eta_0)=\eta_0$.
We choose $k_0$ so that 
\begin{align*}
\sup_{\eta\in\m T\otimes \m H} \big|U(\eta) - U(\pi_{k_0}(\eta))  \big| \leq v=(nh)^{-1/2} \rho,
\end{align*}
which is satisfied if $2^{k_0}=\lfloor \sqrt{n/h}\rfloor$, since 
\begin{align*}
\abs{U(\eta) -  U(\eta')}  \leq \sum_{i=1}^n |W_i| \leq h^{-1} \big[ \norm{g - g'}_{\infty}  + h^{-\kappa} \norm{g}_{\infty} |t - s|^{\kappa}  \big] \leq \sqrt{h/n}\, c\delta(\eta,\eta') \leq v,
\end{align*}
for $\eta'=\pi_{k_0}(\eta)$ by the definition of $\m A_{k_0}$ in their Theorem 5.1 as an $2^{-k_0} v$-net of $\m T\otimes \m H$ relative to the $c\delta$ metric. Under such a choice of $k_0$, the proof of Theorem 5.1 therein leads to the following concentration inequality for any $\alpha> 1/2$ (this truncated version is critical if $\alpha\in(1/2, 1]$, since otherwise the bound of $H_n$ below diverges to $+\infty$ in the original version where $k_0=+\infty$),
\begin{align*}
P \bigg( \sup_{t \in T; \, g \in \m G_n} |U(t, g)| >  H_n+ v\, \sqrt{x}  +  b\, x \bigg) \leq e^{-x}, \quad x > 0, 
\end{align*}
where
\begin{align*}
H_n &\lesssim   v \big\{1 + \sqrt{\log n} + A_n^{1/(2\alpha)} \sum_{k=0}^{k_0} 2^{-k +k/(2\alpha)}\big\}  +  
b \big\{1 + \log n+ A_n^{1/\alpha} \sum_{k=0}^{k_0} 2^{-k +k/\alpha}\big\}\\
& \leq C_1 v\,\big\{1 + \sqrt{\log n} + A_n^{1/(2\alpha)}\big\} + C_2b\, \big\{1 + \log n+ A_n^{1/\alpha} \max\big(1,\, (n/h)^{1/(2\alpha) - 1/2}\big)\big\},
\end{align*}
which proves the desired inequality.

\bibliographystyle{plain}
\bibliography{supbib}

\end{document}